\documentclass[preprint,12pt]{elsarticle}
\usepackage{color}
\usepackage[usenames,dvipsnames,svgnames,table]{xcolor}
\usepackage{amsmath, amsthm, amssymb}

\usepackage{graphicx}
\usepackage{epstopdf}
\def\R{{\mathbb R}}

\newtheorem{problem}{Problem}[section]
\newtheorem{proposition}{Proposition}[section]
\newtheorem{theorem}{Theorem}[section]
\newtheorem{remark}{Remark}[section]

\begin{document}

\begin{frontmatter}

\title{Stability of the kinematically coupled $\beta$-scheme for fluid-structure interaction problems in hemodynamics}

\author[UH]{Sun\v{c}ica \v{C}ani\'{c}\corref{cor1}}
\ead{canic@math.uh.edu}
\author[UH]{Boris Muha}
\ead{borism@math.hr}
\author[UH]{Martina Buka\v{c}}
\ead{martina@math.uh.edu}
\cortext[cor1]{Corresponding author}
\address[UH]{Department of Mathematics,
    University of Houston, Houston, Texas 77204-3476}

\begin{abstract}
It is well-known that classical Dirichlet-Neumann  loosely coupled partitioned schemes for fluid-structure interaction (FSI) problems are unconditionally
unstable for certain combinations of physical and geometric parameters that are relevant in hemodynamics.
It was shown in  \cite{causin2005added} on a simple test problem, that these instabilities are associated 
with the so called  ``added-mass effect''. 
By considering the same test problem as in \cite{causin2005added},
the present work shows that a novel, partitioned, loosely coupled scheme, recently introduced in 
 \cite{MarSun}, called the kinematically coupled $\beta$-scheme, 
 {does not suffer} from the added mass effect for any $\beta \in [0,1]$,
 and is unconditionally stable for all the parameters in the problem.
 Numerical results showing unconditional stability are presented for a full, nonlinearly coupled benchmark FSI problem, first
 considered in~\cite{formaggia2001coupling}. 
\end{abstract}


\end{frontmatter}

\section{Introduction}
Fluid-structure interaction (FSI) problems have important applications in various areas including bio-fluids 
and aero-elasticity. They have been extensively studied from the numerical, as well as analytical point of view
\cite{badia,BarGruLasTuff,bazilevs,bazilevs2,BdV1,
Cervera,ChenShkoller,CSS2,deparis,DEGLT,fernandez,figueroa,gerbeau,heil,Kuk,loon,matthies,BorSun,nobile,peskin1,peskin2,tallec,zhao}. 
A set of popular numerical schemes for FSI in blood flow includes partitioned schemes (loosely or strongly coupled).
Partitioned schemes typically solve an underlying multi-physics problem
by splitting the problem into sub-problems determined by the different physics in the coupled problem.
In particular, in fluid-structure interaction problems the fluid dynamics and structure elastodynamics 
are often solved using separate solvers.
In loosely coupled schemes only one iteration between the fluid and structure sub-problem is performed at each time step,
while in strongly coupled schemes several sub-iterations between the fluid and structure sub-problems need to be performed at each time step
to achieve stability. 

The main advantages of {\sl loosely coupled partitioned schemes} are modularity, simple implementation, and low computational costs. However, in \cite{causin2005added}
it was proved that for certain combinations of physical and geometric parameters (which are realistic in blood flow) ``classical'' 
Dirichlet-Neumann loosely
coupled schemes are unconditionally unstable. In the same paper, the authors showed that this instability is due to the ``added-mass effect''.
Namely, it was shown that a portion of the fluid load to the structure in the coupled FSI problems can be written as an additional inertia term in the structure equation
coming from the fluid mass
(added mass). 
In numerical schemes in which this term appears {\sl explicitly}, which is the case in the classical Dirichlet-Neumann loosely coupled partitioned schemes,
the added mass term acts as a source of instabilities
when the structure is  too "light" to counter-balance the kinetic energy of the "heavy" fluid load.

To get around these difficulties, several different loosely coupled algorithms have been proposed.
The method proposed in~\cite{badia} uses a simple membrane model for the structure 
which can be easily embedded into the fluid problem where it appears as a generalized Robin boundary condition. 
In this way the original problem reduces to a sequence of fluid problems with a generalized Robin boundary condition which can be solved using only the fluid solver.
A similar approach was proposed in~\cite{nobile} where the fluid and structure were split in the classical way, but the fluid and structure sub-problems were linked via novel
transmission (coupling) conditions that improve the convergence rate. 
A different approach to stabilization of loosely coupled (explicit) schemes was proposed in~\cite{burman2009stabilization}
based on the {Nitsche's method}~\cite{hansbo2005nitsche}.
We further mention the scheme proposed in~\cite{badia2009robin} where a Robin-Robin type preconditioner was combined with Krylov iterations for
a solution of an interface system.

For completeness, we also mention several semi-implicit FSI schemes. 
The schemes proposed in~\cite{fernandez2006projection,astorino2009added,astorino2009robin} separate the computation of fluid velocity from the coupled pressure-structure velocity system, 
thereby reducing the computational costs. Similar schemes,  derived from algebraic splitting, were proposed in~\cite{badia2008splitting,quaini2007semi}.
We also mention~\cite{murea2009fast} where an optimization problem is solved at each time-step to achieve continuity of stresses and continuity of velocity at the interface.

Recently, a novel  loosely coupled partitioned scheme, called the ``kinematically coupled $\beta$-scheme'',
was introduced in \cite{MarSun}. 
{Because of its simple implementation, modularity, and good performance,
the kinematically-coupled scheme and its modifications provide an appealing way to study 
multi-physics problems involving FSI. Indeed, this scheme has been used by several groups to study FSI problems
in hemodynamics including poroelastic arterial walls \cite{Martina_Biot},
non-Newtonian fluids \cite{Lukacova},
cardiovascular stents \cite{BorSunStent},  thin structures with longitudinal displacement \cite{CVET,BorSun3D}, FSI with thick structures \cite{thick_structure},
and FSI with multi-layered structure of arterial walls \cite{BorSunMulti,multi-layered}. } 
See also \cite{Fernandez1,Fernandez2,Fernandez3} for a generalization of this scheme called
``the incremental displacement-correction scheme."
The kinematically-coupled $\beta$-scheme successfully deals with problems associated with the added mass effect
in a way different from those reported above. 
It is a modification of the kinematically coupled scheme first introduced in \cite{guidoboni2009stable}.
The parameter $\beta$  was introduced
in \cite{MarSun} to increase the accuracy. 
This parameter distributes the fluid pressure between the fluid and structure sub-problems. 

\if 1 = 0
The first version of the scheme, presented in \cite{MarSun}, concerns FSI problems in which the structure is described by a thin structure model (a thin membrane or shell model; elastic or viscoelastic). Extensions of the scheme have recently been developed by the authors 
to handle FSI problems  with thick structures modeled by the equations of 2D or 3D elasticity \cite{thick_structure}, and to 
handle FSI problems with multiple layered structures (thin-thick-thin-thick\dots; elastic or viscoelastic) \cite{multi-layered}.
Modeling FSI with multi-layered structures is particularly useful in the blood flow application since the walls of major arteries are composed of
several layers: the {\sl tunica intima}, the {\sl tunica media}, and the {\sl tunica adventitia},
 separated by the thin elastic laminae. 
 \fi
 
In \cite{BorSun} the authors used the kinematically-coupled scheme to prove the existence of a weak solution
to a fully nonlinear FSI problem between an incompressible, viscous fluid and a thin structure modeled by either the elastic or viscoelastic
shell equations.
The existence proof is based on constructing approximate solutions using the kinematically-coupled scheme,
and showing that the approximate solutions converge to a weak solution as the time-discretization tends to zero.
This existence result relies on energy estimates, which show that in the kinematically-coupled scheme
the energy of the coupled FSI problem is well-approximated by the discretized problem. 
A consequence of this result is that the kinematically-coupled scheme is stable.
\if 1 = 0The stability of the kinematically coupled scheme is achieved by a clever splitting of the structure problem into 
the inertia part (combined with the viscoelastic part when the (thin) structure is viscoelastic)
and the  elastic part. The structure inertia is then used as a boundary condition in the fluid sub-problem, while the elastodynamics part is solved
separately. The inclusion of the structure inertia (and possibly the viscoelastic part of a thin structure) into the fluid sub-probem was enabled by an implicit
enforcement of the kinematic coupling condition, which enforces the no-slip condition at the fluid-structure interface. {\color{blue} Similar ideas were used by the same authors to prove the existence of a weak solutions to the fluid-multi-layered structure problem in \cite{BorSunMulti}.}
\fi

In contrast with the energy estimates approach showing unconditional stability of the scheme, the present manuscript uses 
a different approach, similar to that of \cite{causin2005added} (which was based on the Dirichlet-to-Neumann mapping applied to a simplified FSI problem), 
which reveals, explicitly, how and why this scheme is unconditionally stable for any $0\le \beta\le 1$. 
\if 1 = 0
The method of proof is different from the one presented in \cite{BorSun}.
It is based on the approach presented in \cite{causin2005added} where it was shown that the classical Dirichlet-Neumann loosely-coupled scheme
is unconditionally unstable. The result in \cite{causin2005added} was proved on the simplest  FSI problem which still retains all the main features
responsible for the stability issues associated with the classical loosely coupled schemes.
The present manuscript considers the same test problem for which it is shown that the kinematically-coupled $\beta$-scheme is
unconditionally stable for all $0\le\beta\le 1$.
\fi

The theoretical results are confirmed by numerical simulations, applied to
a fully nonlinear FSI benchmark problem in hemodynamics, originally presented in \cite{formaggia2001coupling},
and used for testing of various numerical schemes \cite{badia,nobile2001numerical,badia2008splitting,quaini2009algorithms,guidoboni2009stable}.
The hemodynamics and structure parameters in this benchmark problem fall into the critical regime
for which the classical loosely-coupled Dirichlet-Neumann schemes are unstable.
Our results show that the kinematically-coupled $\beta$-scheme is stable for those parameters,
as predicted by the theory. 
We also show that the kinematically-couped $\beta$-scheme 
compares well with the simulations obtained using a monolithic scheme by Badia, Quaini, and Quarteroni \cite{badia2008splitting,quaini2009algorithms}.
Further numerical results obtained using the kinematically-coupled $\beta$-scheme, and a comparison with experimental measurements 
can be found  in  \cite{CVET}.

\if 1 = 0
We conclude this section by mentioning the most recent generalization of the 
kinematically-coupled scheme by Fern{\'a}ndez et al., called
``the incremental displacement-correction scheme'' \cite{Fernandez1,Fernandez2,Fernandez3}.
In \cite{Fernandez3} the authors prove convergence of this numerical scheme by analyzing a simplified
linear model problem. The incremental displacement-correction scheme threats the structure displacement explicitly in the fluid sub-step 
and then corrects it in the structure sub-step. This scheme can also be viewed as a kinematic perturbation of a semi-implicit scheme.
This observation was crucial 
for proving convergence of the scheme \cite{Fernandez3}. 
A different approach, however, needs to be used to prove stability of the 
kinematically coupled $\beta$-scheme due to its particular splitting of the normal stress. 
This gives rise to the various difficulties in estimating the pressure terms which cannot be handled
by using the classical methods presented in \cite{Fernandez3}.



Therefore, the present manuscript provides an original look at the study of stability of the kinematically coupled $\beta$-scheme. 
The method of proof is different from those used in \cite{Fernandez3,BorSun}.
The results of this manuscript provide the information about the stability of the kinematically coupled $\beta$-scheme 
and its relation to the added mass effect which is not captured by the results in \cite{Fernandez3,BorSun}.
This is achieved by an ``explicit calculation'' of the solution to the simplified FSI problem considered here.
The main contribution of this work is in explicitly showing how the implicit enforcement of the kinematic coupling condition
 and the inclusion of the structure inertia into the fluid sub-problem  
 avoid the presence of the added mass effect in the kinematically-coupled $\beta$-scheme.
 \fi

This manuscript is organized as follows.
We begin in Section \ref{FullProblem} by a description of the kinematically-coupled $\beta$-scheme for 
a full FSI problem  in which the structure is modeled
by a cylindrical, linearly {viscoelastic membrane shell model}, also known as the 1D generalized string model.
Then, the simplified FSI problem from \cite{causin2005added} is introduced in Section  \ref{SimplifiedProblem}, 
and a formulation of the kinematically coupled $\beta$-scheme applied to this simplified problem is presented.
Stability analysis of the scheme applied to the simplified problem is presented in Section~\ref{Proof}.
Numerical results are presented in Section \ref{NumericalResults}, where a fully nonlinear 
FSI benchmark problem, discussed in Section \ref{FullProblem}, is solved. 
Stability of the scheme is shown for the parameter values well within the range of instability of Dirichlet-Neumann schemes.
A comparison with the stability analysis associated with classical Dirichlet-Neumann loosely coupled schemes is presented in Section~\ref{Conclusions}.

\section{Problem definition for a fully nonlinear FSI problem}\label{FullProblem}
\begin{figure}[ht]
\centering{
\includegraphics[scale=0.8]{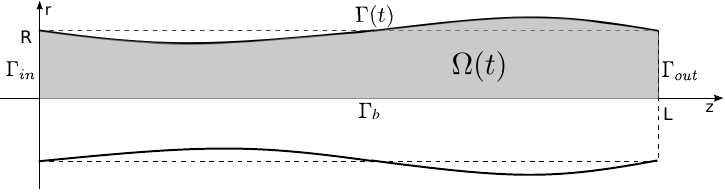}
}
\caption{Deformed domain $\Omega(t).$}
\label{fig:domain}
\end{figure}

We consider the flow of an incompressible, viscous fluid in a channel of radius $R$ and length $L$, see Figure~\ref{fig:domain}.
To fix ideas we will be assuming that the channel is a subset of $\R^2$.
Nothing in the splitting scheme, however, depends on the dimension of the problem,
and the scheme can be applied to 3D problems as well, as was done in, e.g., \cite{BorSun3D}. 
The channel reference domain is denoted by
$\Omega=(0,L)\times (-R,R)$ and the lateral boundary by
$$
\Gamma=\{(z,r)\in \mathbb{R}^2\ | \ 0<z<L,\ r = \pm R\}.
$$
The lateral boundary of the channel is assumed to be deformable, with negligible longitudinal displacement.  
Without loss of generality, we consider only the upper half of the fluid domain with  a symmetry boundary condition at 
the bottom boundary $r=0$. 
The fluid domain, which depends on time, is not known {\sl a priori}, and is denoted by
$$
\Omega(t)=\{(z,r)\in\R^2:z\in (0,L),\ r\in (0,R+\eta(z,t))\},
$$
where $\eta$ denotes the vertical (radial) displacement of the lateral boundary, namely,
$$
\Gamma(t)=\{(z,r)\in\R^2:r=R+\eta(z,t),\ z\in (0,L)\}.
$$
The bottom portion of the boundary will be denoted by $\Gamma_b=(0,L)\times\{0\}$,
with $\Gamma_{in}=\{0\}\times(0,R)$, $\Gamma_{out}=\{L\}\times (0,R)$ denoting the inlet and outlet parts,
respectively.

We consider the flow of an incompressible, viscous fluid driven by the inlet and outlet pressure data $p_{in/out}(t)$,
with the fluid velocity on $\Gamma(t)$ given by $w {\bf e}_r$, and the symmetry boundary condition at the bottom part of the boundary. 
Thus, the {\bf fluid problem} reads:
Find the fluid velocity $\boldsymbol u =(u_z (z,r,t), u_r(z,r,t))$ and pressure $p = p (z,r,t)$ such that 
\begin{equation}
\left \{
\begin{array}{rcll}
\rho_f\big (\displaystyle{\frac{\partial{\bf u}}{\partial t}}+({\bf u}\cdot\nabla){\bf u}\big)&=&\nabla\cdot{\boldsymbol\sigma}  \;  \; &\textrm{in}\; \Omega(t),\ t\in (0,T),  \\
\nabla\cdot{\bf u}&=&0  \;  \; &\textrm{in}\; \Omega(t), \\
u_r=0,\displaystyle{\frac{\partial u_z}{\partial r}}&=&0  \;  \; &\textrm{on}\; \Gamma_b\times (0,T),\\
\boldsymbol\sigma{\bf n}&=&-p_{in/out}(t){\bf n} \;  \; &\textrm{on}\; \Gamma_{in/out}\times (0,T),\\
{\bf u}&=&w{\bf e}_r \;  \; &\textrm{on}\; \Gamma(t),\ t\in (0,T),
\end{array}
\right .
\label{NS}
\end{equation}
where $\rho_f$ is the fluid density, and $\boldsymbol\sigma$ is the fluid stress tensor;
$\boldsymbol\sigma = -p \mathbf{I} + 2 \mu \mathbf{D}(\mathbf{u})$ for a Newtonian fluid, where
$\mu$ is the fluid viscosity and  $\mathbf{D}(\mathbf{u}) = (\nabla \mathbf{u}+(\nabla \mathbf{u})^{\tau})/2$ is the rate-of-strain tensor.

The {\bf structure problem} is defined by solving an elastodynamics problem
for a cylindrical linearly viscoelastic membrane shell, also known as a 1D generalized string model \cite{Quarteroni2000}, capturing only 
vertical (radial) displacement $\eta=\eta(z,t)$ from the reference configuration:
\begin{eqnarray}\label{StructureProblem1}
\rho_s h \frac{\partial^2 \eta}{\partial t^2}+C_0 \eta -C_1 \frac{\partial^2 \eta}{\partial z^2}+D_0\frac{\partial \eta}{\partial t}-
D_1\frac{\partial^3\eta}{\partial t\partial z^2}   =  f, \label{cc2}
\end{eqnarray}
with the boundary conditions $\eta(0)=\eta(L) = 0$, and the initial condition given by 
the zero initial displacement and zero initial structure velocity.
The coefficients $\rho_s$ and $h$ are the structure density and thickness, respectively,
while the constants $C_i, D_i > 0, i = 1,2$, are the elastic and viscoelastic structural coefficients \cite{SunTam,Quarteroni2000}.

The {\bf coupling} between the fluid and structure is defined by the kinematic and dynamics lateral boundary conditions, respectively:
\begin{equation}\label{Coupling1}
\begin{array}{rcl}
\displaystyle{(\frac{\partial \eta}{\partial t}(z,t),0)} &=& \boldsymbol{u}(z,R+\eta(z,),t),\\
f(z,t) &=& -J(z,t)(\boldsymbol\sigma{\bf n})(z,R+\eta(z,t),t)\cdot{{\bf e}_r},
\end{array}
\end{equation}
where $J(z,t)=\displaystyle{-\sqrt{1+\left(\frac{\partial \eta}{\partial z}(z,t)\right)^2}}$
denotes the Jacobian of the transformation from the Eulerian framework used to describe the fluid equations,
to the Lagrangian coordinates used in the structure equations.
Here $f$ in \eqref{Coupling1} is defined by the left hand-side of  \eqref{StructureProblem1}.
The kinematic lateral boundary condition describes continuity of velocities at the fluid-structure interface (the no-slip condition),
while the dynamic lateral boundary condition describes the second Newton's Law of motion on $\Gamma(t)$.

\subsection{The kinematically coupled $\beta$-scheme}
To solve this problem numerically, we consider the kinematically coupled $\beta$-scheme, introduced in \cite{MarSun}.
This scheme is based on the time-discretization via Lie operator splitting  (see e.g., \cite{glowinski2003finite,MarSun}). 

\subsubsection{The Lie scheme}
Let $A$ be an operator from a Hilbert space $H$ into itself, and suppose $\phi_0 \in$ $H$.
Consider the following initial-value problem:
\begin{equation}\label{LieProblem}
\frac{\partial \phi}{\partial t} + A(\phi) = 0, \quad \textrm{in} \ (0,T), \quad {\rm where}\quad  A = \sum\limits_{i=1}^P A_i,
\quad
\phi(0) = \phi_0.
\end{equation}
The Lie scheme consists of splitting the full problem into $P$ sub-problems, each defined by an operator $A_i, i = 1,...,P$.
The original problem is discretized in time with the time step $\Delta t>0$, so that $t^n=n\Delta t$. The Lie splitting scheme consist of solving a series of problems $ \frac{\partial \phi_i}{\partial t} + A_i(\phi_i) = 0$, for $i = 1,...,P$, each defined
over the entire time interval $(t^n, t^{n+1})$, but with the initial data for the $i^{th}$ problem given by the solution of the $(i-1)^{st}$ problem at $t^{n+1}$.
More precisely, set $\phi^0=\phi_0.$ Then, for $n \geq 0$ compute $\phi^{n+1}$ by solving
\begin{equation}
\frac{\partial \phi_i}{\partial t} + A_i(\phi_i) = 0 \quad \textrm{in} \; (t^n, t^{n+1}),
\quad
\phi_i(t^n) = \phi^{n+(i-1)/P}, 
\end{equation}
and then set $\phi^{n+i/P} = \phi_i(t^{n+1}),$ for $i=1, \dots. P.$
This method is first-order accurate in time.
To increase the accuracy in time to second-order,
a symmetrization of the scheme can be performed.

To perform the Lie splitting, problem (\ref{NS})-(\ref{Coupling1}) must first be written as a first-order system in time.
To do this, we utilize the kinematic lateral boundary condition and express
the time-derivative of $\eta$ as the trace of the fluid velocity on $\Gamma$, ${\bf{u}}|_{\Gamma}\cdot {\bf e}_r$. The resulting system 
is then given by the following:
\begin{problem}
Find $({\bf u},p,\eta)$ such that
\end{problem}
\begin{equation*}
\left \{
\begin{array}{rcll}
\rho_f\big (\displaystyle{\frac{\partial{\bf u}}{\partial t}}+({\bf u}\cdot\nabla){\bf u}\big)&=&\nabla\cdot\boldsymbol\sigma  \;  \; &\textrm{in}\; \Omega(t),\ t\in (0,T),  \\
\nabla\cdot{\bf u}&=&0  \;  \; &\textrm{in}\; \Omega(t), \\
u_r=0,&\displaystyle{\frac{\partial u_z}{\partial r}}=0&  \;  \; &\textrm{on}\; \Gamma_b\times (0,T),\\
\boldsymbol\sigma{\bf n}&=&-p_{in/out}(t){\bf n} \;  \; &\textrm{on}\; \Gamma_{in/out}\times (0,T),\\
{\bf u}|_{r = R + \eta}&=&\displaystyle{\frac{\partial\eta}{\partial t}{\bf e}_r} \;  \; &\textrm{on}\; \Gamma,\ t\in (0,T),\\
\displaystyle{\rho_s h \frac{\partial {\bf{u}}}{\partial t}\big|_{\Gamma}\cdot {\bf{e}}_r+C_0 \eta -C_1 \frac{\partial^2 \eta}{\partial z^2}}& &\\
+\displaystyle{D_0 {\bf{u}}|_{\Gamma}\cdot {\bf{e}}_r 
-D_1\frac{\partial^2 ({\boldsymbol{u}}|_{\Gamma})}{\partial z^2}}\cdot{\bf{e}}_r&=&-J\boldsymbol\sigma{\bf n}|_\Gamma\cdot {\bf e}_r\;  \; &\textrm{on}\; \Gamma\times (0,T).
\end{array}
\right .
\label{FSI1}
\end{equation*}
To deal with the motion of the fluid domain, an Arbitrary Lagrangian-Eulerian (ALE) approach is used~\cite{hughes1981lagrangian,donea1983arbitrary,nobile2001numerical}.
We introduce a family of (arbitrary, invertible, smooth) mappings ${\cal{A}}_t$ defined on 
the reference domain $\Omega$  such that for each 
$t \in [0, T]$, ${\cal{A}}_t$ maps the reference domain $\Omega$ onto the current domain $\Omega(t)$.
Now, the ALE time-derivative of function $f$ defined on $\Omega(t)\times (0,T)$ is given by
\begin{equation}
 \frac{\partial f}{\partial t}\bigg|_{\hat{\boldsymbol x}} =  \frac{\partial {f}}{\partial t}+\mathbf{w} \cdot \nabla f,
\end{equation}
where 
$\mathbf{w}$ denotes the domain velocity given by
\begin{equation}
 \mathbf{w}(z,r,t) = \frac{\partial \mathcal{A}_t}{\partial t}(\hat{z},\hat{r}), \ {\rm where}\  (\hat{z},\hat{r}) =   \mathcal{A}_t^{-1}(z,r). \label{w}
\end{equation}
Problem (\ref{NS})-(\ref{Coupling1}) can now be written in the following first-order ALE form
(for more details please see \cite{MarSun}):

\begin{problem}\label{FSI2}
Find $({\bf u},p,\eta)$ such that
\end{problem}
\begin{equation*}
\left \{
\begin{array}{rcll}
\rho_f\big (\displaystyle{\frac{\partial{\bf u}}{\partial t}}\bigg|_{\hat{\boldsymbol x}}+(({\bf u}-\bf w) \cdot\nabla){\bf u}\big)&=&\nabla\cdot\boldsymbol\sigma  \;  \; &\textrm{in}\; \Omega(t),\ t\in (0,T),  \\
\nabla\cdot{\bf u}&=&0  \;  \; &\textrm{in}\; \Omega(t), \\
u_r=0,&\displaystyle{\frac{\partial u_z}{\partial r}}=0&  \;  \; &\textrm{on}\; \Gamma_b\times (0,T),\\
\boldsymbol\sigma{\bf n}&=&-p_{in/out}(t){\bf n} \;  \; &\textrm{on}\; \Gamma_{in/out}\times (0,T),\\
{\bf u}|_{r = R + \eta} &=&\displaystyle{\frac{\partial\eta}{\partial t}{\bf e}_r} \;  \; &\textrm{on}\; \Gamma,\ t\in (0,T),\\
\displaystyle{\rho_s h \frac{\partial ({\bf{u}}|_{\Gamma})}{\partial t}\cdot {\bf{e}}_r+C_0 \eta -C_1 \frac{\partial^2 \eta}{\partial z^2}}& &\\
+\displaystyle{D_0 {\bf{u}}|_{\Gamma} \cdot{\bf{e}}_r
-D_1\frac{\partial^2 ({\bf{u}}|_{\Gamma})}{\partial z^2}}\cdot {\bf{e}}_r&=&-J\boldsymbol\sigma{\bf n}|_\Gamma\cdot {\bf e}_r\;  \; &\textrm{on}\; \Gamma\times (0,T).
\end{array}
\right .
\end{equation*}

The strategy of the kinematically coupled scheme is to split this problem into
a fluid sub-problem and a structure sub-problem in such a way that the inertia of the thin structure is
coupled with the fluid sub-problem via a ``Robin-type''  boundary condition on $\Gamma(t)$. 
When the structure is viscoelastic, the structural viscosity can be treated together with the 
structure inertia as a part of the same boundary condition
(i.e., all the terms in the structure equation
involving the trace of fluid velocity on $\Gamma$ can be used as a boundary condition for the fluid sub-problem).
The elastodynamics of the structure problem is solved separately in the structure sub-problem.

Additionally,  in the $\beta$-scheme a coefficient $\beta \in [0,1]$ is introduced (independent of time), and the fluid stress 
on $\Gamma$ is split into two parts:
$$
\boldsymbol\sigma{\bf n} = \underbrace{\boldsymbol\sigma{\bf n}+\beta{p\mathbf{n}}}_{(I)}
\underbrace{- \beta {p\mathbf{n}}}_{(II)}\quad {\rm on } \quad \Gamma,
$$
where the case $\beta=0$ corresponds to the classical kinematically coupled scheme introduced in \cite{guidoboni2009stable}.
Part I of the fluid stress is used in the fluid sub-problem, while Part II  (the $\beta$-fraction of the pressure) is used to load the
structure in the elastodynamics sub-problem. 
Thus, the dynamic coupling condition on $\Gamma$
\begin{equation*}
\begin{array}{rcl}
\displaystyle{\rho_s h \frac{\partial ({\bf{u}}|_{\Gamma})}{\partial t} \cdot{\bf{e}}_r} &=&{- C_0 \eta +C_1 \displaystyle{\frac{\partial^2 \eta}{\partial z^2}}}
{+\left(\displaystyle{D_0 {\bf {u}}|_{\Gamma} 
-D_1\frac{\partial^2 ({\bf {u}}|_{\Gamma})}{\partial z^2}}-
J \boldsymbol\sigma{\bf n}|_\Gamma-J \beta{p\mathbf{n}|_\Gamma}\right) \cdot{\bf{e}}_r}\\
&& {+ J \beta {p\mathbf{n}|_\Gamma}\cdot {\bf{e}}_r}
\end{array}
\end{equation*}
is split into the fluid part:
\begin{equation}\label{Part1}
 \displaystyle{\rho_s h \frac{\partial ({\bf{u}}|_{\Gamma})}{\partial t}}\cdot{\bf{e}}_r
={ - \displaystyle{D_0 {\bf{u}}|_{\Gamma} \cdot{\bf{e}}_r
+D_1\frac{\partial^2 ({\bf{u}}|_{\Gamma})}{\partial z^2}} \cdot{\bf{e}}_r-
J \boldsymbol\sigma{\bf n}|_\Gamma \cdot {\bf{e}}_r - J \beta{p\mathbf{n}|_\Gamma}\cdot {\bf{e}}_r,}
\end{equation}
and the structure part:
\begin{align*}
 \displaystyle{\rho_s h \frac{\partial ({\bf{u}}|_{\Gamma})}{\partial t} \cdot{\bf{e}}_r}
&={- C_0 \eta +C_1 \frac{\partial^2 \eta}{\partial z^2}}
{+ J \beta {p\mathbf{n}|_\Gamma}\cdot{\bf e}_r,} \  {\rm with}\\
\displaystyle{\rho_s h \frac{\partial ({\bf{u}}|_{\Gamma})}{\partial t}}\cdot{\bf{e}}_r &= \displaystyle{\rho_s h \frac{\partial^2 \eta}{\partial t^2}}.
\end{align*}
This gives rise to the structure problem in terms of $\eta$, defined on $\Gamma$:
\begin{equation}\label{structure_problem}
\displaystyle{\rho_s h \frac{\partial^2 \eta}{\partial t^2}}
{+  C_0 \eta -C_1 \frac{\partial^2 \eta}{\partial z^2}}
=
{J \beta {p\mathbf{n}|_\Gamma}\cdot{\bf e}_r.}
\end{equation}

The splitting scheme is then defined by the following.
\\
\\
\noindent{\bf {The geometry sub-problem:}}\\
\noindent{\bf Step 0.}
The fluid domain $\Omega(t^n)$ and ALE velocity $\mathbf w^{n+1/3}$ are calculated via, e.g.,:
\begin{equation*}
  \mathcal{A}_{t^n}(\hat{z}, \hat{r}) = \left(\hat{z}, \frac{R}{R+\eta^n(\hat{z},t)}\hat{r}\right), \quad \Omega^f(t^{n}) = \mathcal{A}_{t^n}(\hat{\Omega}), 
  \end{equation*}
  and
  \begin{equation*} 
  \mathbf{w}^{n+1/3} = \left(0, \frac{r^n-r^{n-1}}{\Delta t}\right), \label{ale}
\end{equation*}
where $(\hat{z}, \hat{r}) \in \hat{\Omega}, r^{n} \in \Omega(t^n),$ and $r^{n-1} \in \Omega(t^{n-1}).$
\\
\\
\noindent{\bf {The fluid sub-problem (Steps 1 and 2):}}\\
\noindent{\bf Step 1.}
A time-dependent Stokes problem is solved on the fixed fluid domain $\Omega(t^n)$ obtained in {Step 0}.
The boundary condition at the lateral boundary  is given by \eqref{Part1}, 
with $\beta J {p\mathbf{n}}$ taken explicitly from the previous time step. 
The displacement of the structure stays intact. The problem reads:\\
{\sl
Given $p^n$, $\eta^n$, and $\mathbf{u}^n$ from the previous time step, find ${\bf u}^{n+1/3}$, $p^{n+1}$, $\eta^{n+1/3}$ such that:
\begin{equation}\label{step1}
 \left\{\begin{array}{l@{\ }}
\displaystyle{ \rho_f \frac{ \mathbf{u}^{n+1/3}-\mathbf{u}^n}{\Delta t}} =\nabla \cdot \boldsymbol{\sigma}^{n+1/3},   \quad \nabla \cdot \mathbf{u}^{n+1/3}=0 \quad \textrm{in} \; \Omega(t^n) \\ \\
 \displaystyle{\rho_s h \frac{{\bf{u}}^{n+1/3}-{\bf{u}}^n}{\Delta t}|_{\Gamma} \cdot {\bf{e}}_r
 +D_0 {\bf{u}}^{n+1/3}|_{\Gamma}\cdot {\bf{e}}_r -D_1 \frac{\partial^2{\bf{u}}^{n+1/3}}{ \partial z^2}|_{\Gamma}} \cdot {\bf{e}}_r\\
+ \sqrt{(1+\left(\frac{\partial \eta^n}{\partial z}\right)^2} (\boldsymbol \sigma^{n+1/3} \boldsymbol {n})\cdot \mathbf{e}_r   
= -  {\beta  \sqrt{(1+\left(\frac{\partial \eta^n}{\partial z}\right)^2} p^n {\bf n}\cdot\mathbf{e}_r} 
 \quad \textrm{on} \; \Gamma, \\ \\
 {\bf{u}}^{n+1/3} \cdot {\bf{e}_z} = 0   \quad \textrm{on} \; \Gamma, \\ \\
\displaystyle{\eta^{n+1/3}  = \eta^n }\quad \textrm{in} \; \Gamma,
  \end{array} \right.  \nonumber
\end{equation}
with the following boundary conditions: 
\begin{equation*}
  \frac{\partial u_z^{n+1/3}}{\partial r}(z,0,t) =  u_r^{n+1/3}(z,0,t) = 0 \quad \textrm{on} \; \Gamma_b,\ u_z^{n+1/3}=0\ {\rm on}\ \Gamma(t^n)
\end{equation*}
\begin{equation*}
  \mathbf{u}^{n+1/3}(0,R,t) = \mathbf{u}^{n+1/3}(L,R,t) = 0, 
\end{equation*}
\begin{equation*}
 \boldsymbol\sigma^{n+1/3} \mathbf{n}|_{in} = -p_{in}(t^{n+1})\mathbf{n}|_{in}\  {\rm on}\ \Gamma_{\rm in}, \; \; \boldsymbol\sigma^{n+1/3} \mathbf{n}|_{out} = -p_{out}(t^{n+1})\mathbf{n}|_{out}  \ {\rm on}\ 
 \Gamma_{\rm out}.
\end{equation*}
}
\\
\textbf{Step 2.} Solve the fluid and ALE advection sub-problem defined on $\Omega(t^n)$,
with the ALE velocity ${\bf{w}}^{n+1/3}$ defined by \eqref{ale} based on the domain $\Omega(t^n)$ and the corresponding ALE mapping.
The problem reads:\\
{\sl Find $\mathbf{u}^{n+2/3}$ and ${\eta}^{n+2/3}$
such that 
\begin{equation*}
 \left\{\begin{array}{l@{\ }} 
 \displaystyle{\frac{\mathbf{u}^{n+2/3}-\mathbf{u}^{n+1/3}}{\Delta t}} + (\mathbf{u}^{n+1/3}-\mathbf{w}^{n+1/3}) \cdot \nabla \mathbf{u}^{n+2/3}= 0,   \quad \textrm{in} \; \Omega(t^n)  \\ \\
 \displaystyle{{\eta}^{n+2/3} = \eta^{n+1/3}} \quad \textrm{on} \; \Gamma,\\ \\
\displaystyle{{\bf{u}}^{n+2/3}|_{\Gamma}={\bf u}^{n+1/3}|_{\Gamma} } \quad \textrm{on} \; \Gamma, \end{array} \right. 
\end{equation*}
with the following boundary conditions: 
$$\mathbf{u}^{n+2/3}=\mathbf{u}^{n+1/3} \  \; \textrm{on} \; \Gamma_{-}^{n+1/3}(t^n), \; \textrm{where}$$ 
$$\Gamma_{-}^{n+1/3}(t^n) = \{\mathbf{x} \in \mathbb{R}^2 | \mathbf{x} \in \partial \Omega(t^n), (\mathbf{u}^{n+1/3}-\mathbf{w}^{n+1/3})\cdot \mathbf{n} <0 \}.$$}

\noindent{\bf{The structure (elastodynamics) sub-problem (Step 3):}}
\\
\noindent\textbf{Step 3.}  Solve the elastodynamics problem \eqref{structure_problem} for the location of the deformable boundary
by involving the elastic part of the structure which is loaded by Part II of the normal fluid stress.
Additionally, the fluid and structure communicate via the kinematic lateral boundary condition which gives
the velocity of the structure in terms of the trace of the fluid velocity, taken initially to be the value from the previous step. To discretize the elastodynamics problem in time, we use the $\theta$--scheme (see~\cite{glowinski2003finite} for details).  
The problem reads:\\
{\sl
Given $p^{n+1}$ computed in Step 1, and $\eta^n$ obtained at the previous time step, find ${\mathbf{u}^{n+1}}$ and ${\eta^{n+1}}$ such that  
\begin{equation*}
 \left\{\begin{array}{l@{\ }} 
\displaystyle{\mathbf{u}^{n+1} = \mathbf{u}^n},   \quad \textrm{in} \; \Omega(t^n) \\ \\
 \displaystyle{ {\bf{u}}^{n+1}|_{{\Gamma(t^n)}} \cdot {\bf{e}}_r=\frac{\eta^{n+1}-\eta^{n-1}}{2\Delta t}} \quad \textrm{on} \; \Gamma,\\ \\
\displaystyle{\rho_s h \frac{\eta^{n+1}-2\eta^{n}+\eta^{n-1}}{\Delta t^2}+C_0 (\theta\eta^{n+1}+(1-2\theta)\eta^{n}+\theta\eta^{n-1})} \\
\displaystyle{-C_1\left(\theta\frac{\partial^2 \eta^{n+1}}{\partial z^2}+(1-2\theta)\frac{\partial^2 \eta^{n}}{\partial z^2}+\theta\frac{\partial^2 \eta^{n-1}}{\partial z^2}\right)
=  \beta \sqrt{(1+\left(\frac{\partial \eta^n}{\partial z}\right)^2} p^{n+1} \mathbf{n}\cdot\mathbf{e}_r}  \quad \textrm{on} \; \Gamma, \end{array} \right.  
\end{equation*}
with the boundary conditions: 
\begin{equation*}
  \eta^{n+1}|_{z=0,L}=0.
\end{equation*}
}
Set $t^n=t^{n+1}$, and return to Step 0.

\begin{figure}[ht]
 \centering{
 \includegraphics[scale=0.4]{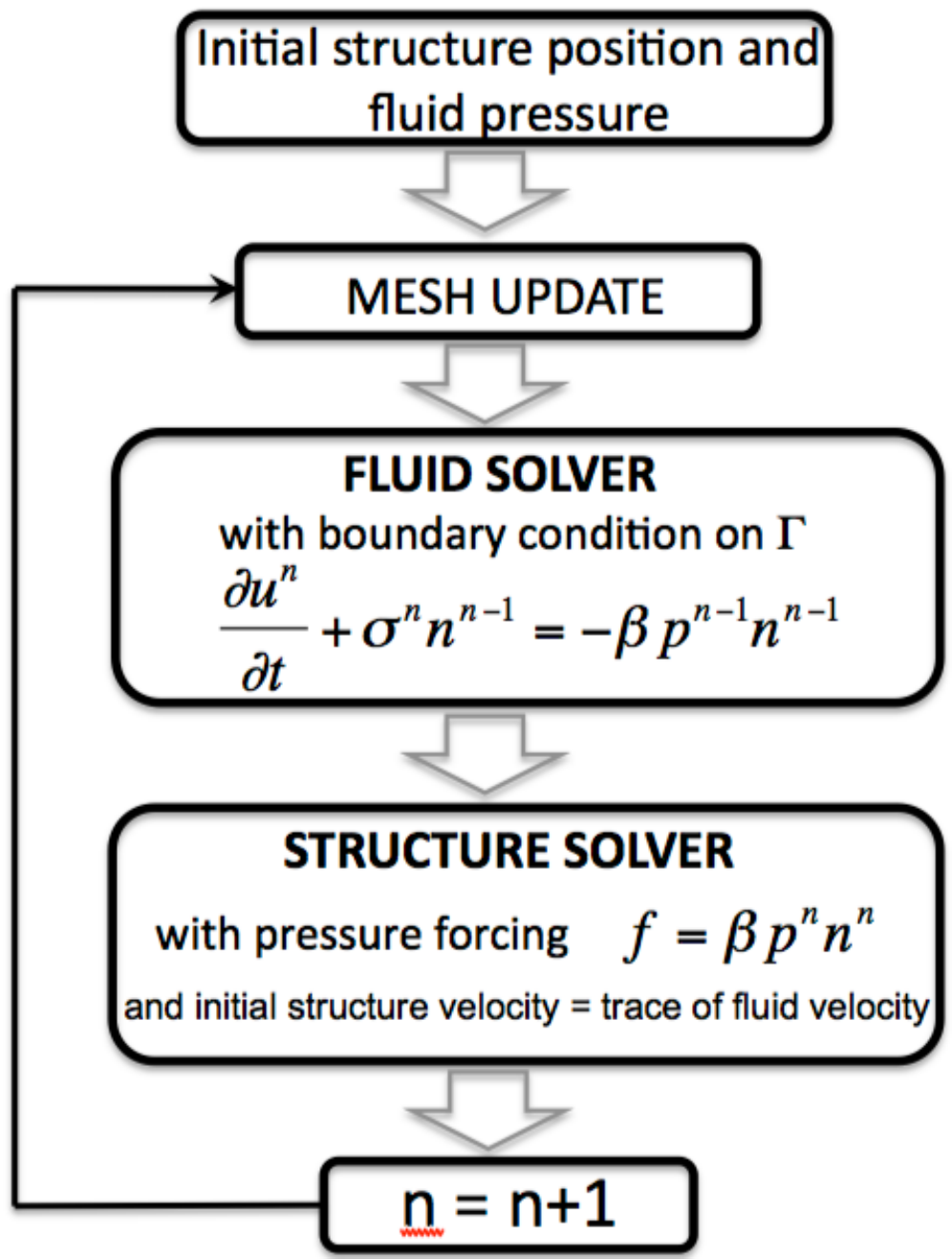}
 }\label{kk_beta_scheme}
 \caption{A block diagram showing the main steps of the kinematically coupled $\beta$-scheme.}
 \end{figure}

A block-diagram showing the main steps of the scheme is given in Figure~\ref{kk_beta_scheme}.
More details about the scheme can be found in \cite{MarSun}.

\section{The simplified test problem}\label{SimplifiedProblem}
We show that the 
 kinematically coupled $\beta$-scheme is unconditionally stable 
when applied to the simplified problem considered in \cite{causin2005added} as the simplest problem which captures
 the main features related to the instabilities in classical loosely-coupled schemes caused by the ``added mass effect''.
The simplified problem consists of solving the time-dependent Stokes equations for an incompressible, inviscid fluid
in a 2D channel with deformable walls, and with the elastodynamics equations given by the 1D generalized string model.
Moreover, it is assumed that the displacement of the deformable wall is small enough so that
it can be neglected in the fluid flow problem. In this case the geometry of the fluid domain is fixed while the small
deformation of the boundary is calculated using the elastodynamics equations, which are coupled with the fluid flow
through the kinematic and dynamic coupling conditions.
Thus, the problem is defined on the reference fluid domain $\Omega$ with the lateral boundary $\Gamma$.

We will be assuming that the inlet and outlet pressure data, which are functions of time, are uniformly bounded. This is a reasonable assumption 
for the blood flow application, since the pressure is typically a periodic function of time, with a bounded amplitude.

The {\bf fluid problem} reads:
Find the fluid velocity $\boldsymbol u = \boldsymbol u (z, r, t)$ and pressure $p = p (z,r,t)$ such that 
\begin{equation}\label{FluidProblem}
 \left\{\begin{array}{ r c l l}
\displaystyle{\rho_f  \frac{\partial \boldsymbol u}{\partial t}}  + \nabla p &=& 0  \;  \; &\textrm{in}\; \Omega \times(0,T),  \\  
\nabla \cdot \boldsymbol {u} &=& 0   \;  \; &\textrm{in}\; \Omega \times(0,T), \\
\boldsymbol u \cdot  \boldsymbol n &=& 0  \;  \; &\textrm{on}\; \Gamma_b \times(0,T), \\
p &=& p_{in/out} (t)  \;  \; &\textrm{on}\; \Gamma_{in/out} \times(0,T), \\
\boldsymbol u \cdot {\bf{e}}_r &=& w    \;  \; &\textrm{on}\; \Gamma \times(0,T),
 \end{array} \right.  
\end{equation}
with the initial velocity and pressure equal to zero.

The {\bf structure problem} is defined by solving for $\eta=\eta(z,t)$ the 1D generalized string model:
\begin{eqnarray}\label{StructureProblem}
\rho_s h \frac{\partial^2 \eta}{\partial t^2}+C_0 \eta -C_1 \frac{\partial^2 \eta}{\partial z^2}   =  f, \label{cc2}
\end{eqnarray}
with boundary conditions $\eta(0)=\eta(L) = 0$, and initial conditions given by the
zero initial displacement and zero initial structure velocity.

\if 1 = 0
\begin{remark}
Notice that the simplified problem does not have the viscoelastic terms. The addition of the viscoelastic terms only contributes
to the stability of the scheme. 
{For more details about the analysis of the fully nonlinear FSI problem with and without the viscoelastic terms see \cite{BorSun}.} 
\end{remark}
\fi

The {\bf coupling} between the fluid and structure is defined by the kinematic and dynamic lateral boundary conditions,
which, in this simplified problem, read:
\begin{equation}\label{Coupling}
\begin{array}{rcl}
\displaystyle{\big(\frac{\partial \eta}{\partial t},0\big) }&=& \boldsymbol{u}|_{\Gamma},\\
\displaystyle{\rho_s h \frac{\partial^2 \eta}{\partial t^2}+C_0 \eta -C_1 \frac{\partial^2 \eta}{\partial z^2}}   &=& p|_{\Gamma}.
\end{array}
\end{equation}

The coupled fluid-structure interaction problem can be written as follows:
\begin{problem}\label{PS1}
\begin{equation}
 \left\{\begin{array}{ r c l l}
\displaystyle{\rho_f  \frac{\partial \boldsymbol u}{\partial t}}  + \nabla p &=& 0  \;  \; &\textrm{in}\; \Omega \times(0,T),  \\  
\nabla \cdot \boldsymbol {u} &=& 0   \;  \; &\textrm{in}\; \Omega \times(0,T), \\
\boldsymbol u \cdot  \boldsymbol n &=& 0  \;  \; &\textrm{on}\; \Gamma_b \times(0,T), \\
p &=& p_{in/out} (t)  \;  \; &\textrm{on}\; \Gamma_{in/out} \times(0,T), \\
\boldsymbol u \cdot   \bf{e}_r &=& \displaystyle{\frac{\partial \eta}{\partial t}}   \;  \; &\textrm{on}\; \Gamma \times(0,T),\\
\displaystyle{\rho_s h \frac{\partial^2 \eta}{\partial t^2}+C_0 \eta -C_1 \frac{\partial^2 \eta}{\partial z^2}} &=& p\;  \; &\textrm{on}\; \Gamma \times(0,T).
 \end{array} \right.  
\label{StabilityStokes}
\end{equation}
\end{problem}
To perform the Lie splitting, this system is written as a first-order system in time by using the kinematic coupling condition to replace 
$\displaystyle{\rho_s h \frac{\partial^2 \eta}{\partial t^2}}$ by $\displaystyle{\rho_s h \frac{\partial \boldsymbol u}{\partial t}}|_{\Gamma} \cdot {\bf{e}}_r$.
As before, this point is crucial in order to perform the splitting which gives rise to a stable scheme. Namely, 
in the fluid sub-problem, which we write next, the structure inertia will be taken into account implicitly in the boundary condition on $\Gamma$,
which is a crucial ingredient for the stability of the $\beta$-scheme. {\bf This is a way how the fluid is coupled to the structure in this splitting.}

To simplify notation, 
as in \cite{causin2005added}, we introduce a linear, symmetric, positive definite operator ${\mathcal{L}}$ defined by the elasticity 
tensor associated with the structure problem. Namely, for $\eta,\xi \in H^1_0(0,L)$, we define
\begin{equation}
<\mathcal L \eta,\xi > := a_S(\eta,\xi),
\end{equation}
where $a_S(\eta,\xi)$ is the inner product on $H^1(0,L)$ defined by
$$
a_S(\eta,\xi) := \int_0^L C_0 \eta\xi dx + \int_0^L C_1 \frac{\partial\eta}{\partial x}\frac{\partial\xi}{\partial x} dx, \quad \forall \eta,\xi \in H^1(0,L).
$$
With this notation, the structure equation in Step 2 can be written as
\begin{equation*}
\rho_s h \frac{\partial^2 \eta}{\partial t^2} + \mathcal L \eta =  p|_{\Gamma}.
\end{equation*}
Notice that we are free to choose for the operator $\mathcal L$ any positive-definite operator associated with
elastodynamics of thin structures, such as, for example, 
 the full Koiter shell model operator, described in \cite{ SunTam}. 

Before we write the main splitting steps we also notice that the fluid stress,
which in this simplified example corresponds to the pressure, will be split into two parts: Part I which is given by $p - \beta p$, and Part 2 which is 
given by $\beta p$, so that $p = (p - \beta p) + \beta p$=Part 1 + Part 2.

Thus, the main steps of the Lie splitting for the simplified FSI problem \eqref{StabilityStokes} are:\\
\\
\noindent{\bf{The fluid sub-problem:}}
\vskip 0.1in
\noindent
{\bf Step 1.} The Stokes problem is solved on a fixed fluid domain, with a boundary condition
which couples the structure inertia with Part I of the fluid stress. Here, the portion $\beta p$ of the stress is taken explicitly,
while the rest is taken implicitly. The displacement of the structure stays intact. The problem reads:
{\sl 
Given $p^n$ from the previous time step,
find $\mathbf{u}^{n+1/2},  \eta^{n+1/2}$ and $p^{n+1}$ such that:
\begin{equation}
 \left\{\begin{array}{ r c l l}
\displaystyle{\rho_f  \frac{\mathbf u^{n+1/2}-\mathbf{u}^n}{\Delta t}}  + \nabla p^{n+1} &=& 0  \;  \; &\textrm{in}\; \Omega,  \\  
\nabla \cdot \mathbf {u}^{n+1/2} &=& 0   \;  \; &\textrm{in}\; \Omega, \\
\mathbf u^{n+1/2} \cdot  \boldsymbol n &=& 0  \;  \; &\textrm{on}\; \Gamma_b, \\
p^{n+1} &=& p_{in/out} (t)  \;  \; &\textrm{on}\; \Gamma_{in/out}, \\
\displaystyle{\rho_s h \frac{ {\mathbf{u}^{n+1/2}} - \mathbf{u}^n}{\Delta t} }\big|_{\Gamma}-p^{n+1} &=& - \beta p^n. \; &\textrm{on}\; \Gamma,
 \end{array} \right.  
 \end{equation}
 with the {\underline{initial data for the fluid velocity}} on $\Gamma$, namely, $u^n|_\Gamma$, given by the structure velocity calculated from the previous time-step.
The structure displacement stays intact, namely
\begin{equation}
\eta^{n+1/2} = \eta^n \quad {\rm on}\quad \Gamma.
\end{equation}
}
\vskip 0.1in
\noindent{\bf{The structure sub-problem:}}
\vskip 0.1in
\noindent
{\bf Step 2.} The structure problem is solved with the pressure load $\beta p^{n+1}$ just calculated from the fluid sub-problem.
The problem reads:
{\sl 
Find ${\mathbf{u}^{n+1}}$ and $\eta^{n+1}$ such that the following holds:
\begin{equation*}
 \left\{\begin{array}{l@{\ }} 
 \displaystyle{\rho_s h \frac{\eta^{n+1}-2\eta^{n}+\eta^{n-1}}{\Delta t^2}+\mathcal{L} (\theta\eta^{n+1}+(1-2\theta)\eta^{n}+\theta\eta^{n-1}) 
=  \beta p^{n+1}|_{\Gamma}}  \quad \textrm{on} \; \Gamma,\\
 \displaystyle{ {\bf{u}}^{n+1}|_{{\Gamma(t^n)}} \cdot {\bf{e}}_r=\frac{\eta^{n+1}-\eta^{n-1}}{2\Delta t}} \quad \textrm{on} \; \Gamma,
  \end{array} \right.   
\end{equation*}
with the {\underline{initial data for the structure velocity}} given by the trace of the fluid velocity on $\Gamma$, obtained from the just-calculated
fluid sub-problem.
The fluid velocity  in $\Omega$ remains intact, namely:
$$
 \mathbf {u}^{n+1}= \mathbf{u}^{n+1/2}   \quad \textrm{in} \; {\Omega}.
$$
}

The coupling of the two sub-problems through the initial data is the classical coupling provided by 
the Lie splitting scheme. 
\vskip 0.1in
\noindent
{\bf The pressure formulation:}
\vskip 0.1in
\noindent
In this simplified model, the problem is Step 1 can be entirely formulated in terms of the pressure and the normal trace of  the fluid velocity on $\Gamma$.
Namely, by taking the divergence free condition in the differentiated first (momentum) equation,
the problem can be written in terms of the Laplace operator for the pressure.
Similarly, we can re-write the boundary conditions on $\Gamma$ in terms of the pressure and the normal trace of the velocity on $\Gamma$. 

To simplify notation we will be using $\boldsymbol{n}$ to denote the normal on $\Gamma$ (equal to ${\bf{e}}_r$),
and $u_r|_{\Gamma}$ to denote the normal trace of the fluid velocity on $\Gamma$:
\begin{equation*}
\boldsymbol{n} := {\bf{e}}_r\ {\rm on}\ \Gamma\quad  {\rm and}\quad  u_r|_{\Gamma} := {\bf{u}}|_{\Gamma} \cdot {\bf{e}}_r.
\end{equation*}
With this, the problem is Step 1 can be stated as:\\
\break
\centerline{\bf{The fluid sub-problem in terms of $p$ and $u_r|_{\Gamma} $}}
\vskip 0.05in
\noindent
{\bf Step 1'.} {\sl Given $u_r^{n}|_{\Gamma}$, $p^n, \eta^n$, 
find $p^{n+1/2}$,  $u_r^{n+1/2}|_{\Gamma}$, $\eta^{n+1/2}$, such that for $t\in (t^n, t^{n+1})$:
\begin{equation}\label{Step1}
 \left\{\begin{array}{rcll}
 \displaystyle{-\triangle p^{n+1/2} } &=& 0,  \quad & \textrm{in} \; \Omega, \\ 
  \displaystyle{p^{n+1/2}   }&=& p_{in/out}(t)  \quad &\textrm{on} \; \Gamma_{in/out}, \\
 \displaystyle{\frac{\partial p^{n+1/2}}{\partial \boldsymbol n}}  &=& 0 \quad &\textrm{on} \; \Gamma_{b}, \\ 
\displaystyle{p^{n+1/2}  + \frac{\rho_s h}{\rho_f} \frac{ \partial p^{n+1/2}}{\partial \boldsymbol n} }&=& \beta p^n    \quad &\textrm{on}\; \Gamma, \\ 
-----&--&-----&---\\
\displaystyle{\frac{ u_r^{n+1/2}|_{\Gamma}-u_r^{n}|_{\Gamma}}{\Delta t}} &=& - \displaystyle{\frac{1}{\rho_f}\frac{ \partial p^{n+1/2}|_{\Gamma}}{\partial \boldsymbol n}} \quad &\textrm{on}\; \Gamma,\\
\displaystyle{\eta^{n+1/2}}&=& \displaystyle{\eta^n} \quad &\textrm{on}\; \Gamma,
\end{array} \right. 
\end{equation}
}with the {\underline{initial data for the fluid velocity}} on $\Gamma$ given by the structure velocity calculated from the previous time-step.

The first four equations in Step 1
define a time-dependent Robin problem for the pressure
(supplemented with the initial data for the pressure). 
The time-dependence enters through the inlet and outlet pressure data, which are functions of time.
Notice that this problem already incorporates a portion of the coupling of the underlying FSI problem.
This problem has a unique solution in the space $H^1(\Omega)$ provided that the data
$p_{in/out} \in H^{1/2}(\Gamma_{in/out})$, and $p^n\in H^{-1/2}(\Gamma)$.

One can also show (see e.g., \cite{causin2005added}) that the problem for the structure, defined in Step 2, has a unique solution 
$\eta \in H^1_0(0,L)$ provided that the initial data are such that $\eta|_{t = 0} \in H^1_0(0,L)$, and $\partial\eta/\partial t|_{t = 0} \in H^1_0(0,L)$,
with $p^{n+1}|_{\Gamma} \in H^{1/2}(\Gamma)$.

\begin{remark}
We emphasize that the fluid sub-problem in Step 1' is coupled to the structure sub-problem via the following two equations:
\begin{itemize}
\item The Robin boundary condition, given by the 4th equation in \eqref{Step1}. 
This term accounts for the structure inertia,
in an implicit way, via the kinematic coupling condition, and a portion of the dynamic coupling condition on $\Gamma$. 
This is exhibited, in particular, through the term 
$\frac{\rho_s h}{\rho_f} \frac{ \partial p^{n+1/2}}{\partial \boldsymbol n}$.
In this term, $\rho_s h$ is the structure mass, and $\frac{ \partial p^{n+1/2}}{\partial \boldsymbol n}$,
which was calculated from 
$\frac{ \partial p^{n+1/2}}{\partial \boldsymbol n}= -\rho_f \frac{\partial^2 \eta^{n+1/2}}{\partial t^2}$, accounts for the structure inertia 
\underline{implicitly}. If the structure inertia were taken into account explicitly, as is the case in Dirichlet-Neumann schemes, one would have the
Neumann boundary condition $\frac{ \partial p^{n+1/2}}{\partial \boldsymbol n}= -\rho_f \frac{\partial^2 \eta^{n}}{\partial t^2}$
for the pressure holding on $\Gamma$, given {\underline{explicitly}} in terms of the structure velocity from the previous time-step,
instead of the Robin condition. 
Our Robin boundary condition comes from: (1) using a portion of the dynamic coupling condition as dictated by the Lie splitting,
through which the structure inertia $\rho_s h \frac{\partial^2 \eta}{\partial t^2}$ is included in the fluid sub-problem, 
and (2)  using the kinematic coupling condition $\frac{\partial^2 \eta}{\partial t^2} = {\bf u}|_{\Gamma}\cdot {\bf e}_r$
to replace the structure velocity by the trace on the fluid velocity on $\Gamma$ in an {\underline {implicit}} way.
By writing the problem in terms of the pressure, this boils down to using the condition 
$\frac{ \partial p^{n+1/2}}{\partial \boldsymbol n}= -\rho_f \frac{\partial^2 \eta^{n+1/2}}{\partial t^2}$ in an implicit way,
and not explicitly as in the Dirichlet-Neumann schemes.
The implicit use of this condition leads to the fact that our fluid sub-problem, which includes structure inertia implicitly, can be written entirely 
in terms of the pressure. 
Thus, in this fluid sub-problem, the fluid and a portion of the structure equation are implicitly coupled
via the implicit enforcement of the kinematic coupling condition, which leads to a fluid sub-problem for the pressure with a Robin-type boundary 
condition, which includes the structure inertia through the term multiplying $\rho_s h/\rho_f$.
\item The initial condition for the trace of the fluid velocity on $\Gamma$, which is given by the structure velocity calculated from the previous 
time-step.
\end{itemize}

Likewise, the structure problem in Step 2 is coupled to the fluid problem given in Step 1' 
via:
\begin{itemize}
\item the pressure loading term $\beta p^{n+1}|_{\Gamma}$, and 
\item through the initial data for the structure velocity, which is given
by the trace of the fluid velocity from the just-calculated fluid sub-problem. The latter reflects the ``standard'' coupling
provided by the Lie splitting scheme. Thus, even for $\beta = 0$, the structure ``feels'' the fluid via the initial condition for the 
structure velocity.
\end{itemize}
\end{remark}


\section{Stability analysis}\label{Proof}

To study stability of the kinematically-coupled $\beta$-scheme we introduce the following operator:
let $\mathcal P : H^{-1/2}(\Gamma) \to Q$,
where 
\begin{equation}\label{Q}
Q = \{ q \in H^1(\Omega) \ | \ q|_{\Gamma_{in/out}} = 0\},
\end{equation}
such that
\begin{equation}\label{P}
 \left\{\begin{array}{rcll}
 \displaystyle{-\triangle \mathcal {P}w }&=& 0,   \quad  &\textrm{in} \; \Omega, \\ 
  \displaystyle{\mathcal {P}w   }&=& 0  \quad &\textrm{on} \; \Gamma_{in/out}, \\
 \displaystyle{\frac{\partial \mathcal {P}w}{\partial \boldsymbol n}}  &=& 0 \quad &\textrm{on} \; \Gamma_{b}, \\[0.3cm] 
\displaystyle{\mathcal {P}w + \frac{\rho_s h}{\rho_f} \frac{ \partial \mathcal {P}w}{\partial \boldsymbol n} } &=& w    \quad&\textrm{on}\; \Gamma.
\end{array} \right. 
\end{equation}
Operator $\mathcal P$ associates to every $w\in H^{-1/2}(\Gamma)$ the solution of the pressure problem in Step 1',
with the homogeneous inlet and outlet data $p_{in/out} = 0$.

We are interested in the trace on $\Gamma$ of this pressure solution. For this purpose, 
introduce the operator $\mathcal S : H^{-1/2}(\Gamma) \to H^{1/2}(\Gamma)$ by
\begin{equation}
\mathcal S w = \mathcal P w|_{\Gamma}.
\end{equation}
One can prove  the following classical results (see e.g. \cite{causin2005added,QuarteroniValli}):

\begin{proposition}
Operator $\mathcal S$ satisfies the following properties:
\begin{enumerate}
\item  $\mathcal S:  H^{-1/2}(\Gamma) \to H^{1/2}(\Gamma)$ is continuous.
\item  $\mathcal S$ is compact, self-adjoint, and positive on $L^2(\Gamma)$.
\end{enumerate}
\end{proposition}
\noindent
{\bf Proof.}
Since this result is standard, we just outline the proof.
We start by writing the weak formulation of problem (\ref{P}) which reads: find ${\cal P}w\in Q$ such that
\begin{equation}
\int_{\Omega}\nabla{\cal P}w:\nabla\phi+\frac{\rho_f}{\rho_s h}\int_{\Gamma}{\cal P}w\phi=\frac{\rho_f}{\rho_s h}\langle w,\phi\rangle,\quad \forall \phi\in Q,
\label{PWF}
\end{equation}
where $\langle w,\phi\rangle$ is a duality pairing between $H^{-1/2}(\Gamma)$ and $H^{1/2}(\Gamma)$.
The first assertion (continuity) follows directly from here.

To show the second statement of the proposition, we consider  $w,v\in L^2(\Gamma)$ and notice:
\begin{equation}
\int_{\Gamma}{\cal S}wv=\int_{\Gamma}{\cal P}wv=\frac{\rho_s h}{\rho_f}\int_{\Omega}\nabla{\cal P}v:\nabla{\cal P}w+\int_{\Gamma}{\cal P}v{\cal P}w.
\label{SymS}
\end{equation}
Here, the first equality follows from the definition of operator ${\cal S}$. For the second equality we notice that ${\cal P}w\in Q$ is an admissible test function
and, therefore, the second equality is just formula (\ref{PWF}) with the test function $\phi={\cal P}w$. 
We now consider ${\cal P}$ as an operator on $L^2(\Gamma)$, i.e., ${\cal S}:L^2(\Gamma)\rightarrow L^2(\Gamma)$.
For every $w\in L^2(\Gamma)$, we have ${\cal P}w\in Q$. The trace theorem implies that ${\cal S}w={\cal P}w|_{\Gamma}\in H^{1/2}(\Gamma)$,
and so 
we have Im$({\cal S})\subseteq H^{1/2}(\Gamma)$. Since $\Gamma$ is a bounded set, embedding $H^{1/2}(\Gamma)\hookrightarrow L^2(\Gamma)$ is compact,
and thus ${\cal S}$ is a compact operator. Formula (\ref{SymS}) implies that operator ${\cal S}$ is positive and self-adjoint,
which completes the proof.
\qed

To study the solution of the corresponding non-homogeneous problem with $p=p_{in/out}(t)$ on $\Gamma_{in/out}$,
we introduce an arbitrary continuous extension operator $E_F : H^{1/2}(\partial \Omega \backslash \Gamma) \rightarrow H^1(\Omega)$ so that $E_F q|_{\partial \Omega \backslash \Gamma} = q$ and $||E_F q||_{H^1(\Omega)} \le C ||q||_{H^{1/2}(\partial \Omega \backslash \Gamma)}. $ Let $p^* \in C(0, \infty; H^1(\Omega))$ be the solution to
\begin{equation}
 \left\{\begin{array}{rcll}
 \displaystyle{-\triangle p^*}  &=& \triangle E_F \bar{p},  \quad & \textrm{in} \; \Omega, \\ 
  \displaystyle{p^*  } &=& 0 \quad &\textrm{on} \; \Gamma_{in/out}, \\ 
 \displaystyle{\frac{\partial p^*}{\partial \boldsymbol n}}  &=& -\displaystyle{\frac{\partial E_F \bar{p}}{\partial \boldsymbol n}}   \quad &\textrm{on} \; \Gamma_{b}, \\ 
\displaystyle{p^*  + \frac{\rho_s h}{\rho_f} \frac{ \partial p^*}{\partial \boldsymbol n}} &=& -\displaystyle{\frac{\partial E_F \bar{p}}{\partial \boldsymbol n}}   \quad &\textrm{on}\; \Gamma,
\end{array} \right.  \nonumber
\end{equation}
where $\bar{p} = p_{in/out}$ on $\Gamma_{in/out}$.
Now, the solution to the pressure problem in Step 1'  is given by
\begin{equation*}
p = p^* +E_F \bar{p} + \mathcal P (\beta p^n).
\end{equation*}
By denoting
\begin{equation*}
p_{ext} = p^*|_{\Gamma}+E_F \bar{p}|_{\Gamma}
\end{equation*}
{we can write the trace of the pressure solution in Step 1' on $\Gamma$ as}
\begin{equation}
 p|_{\Gamma} = p_{ext} + \mathcal S (\beta p^n).
\label{pext}
\end{equation}
This trace of the pressure is used to load the equation for the structure in Step 2. 
Thus, the structure problem in Step 2 can now be written as:
\textit{Find $\eta$ such that}
\begin{equation} 
\rho_s h \frac{\partial^2 \eta}{\partial t^2} + \mathcal{L} \eta =  \beta (p_{ext}^{n+1} + \mathcal S (\beta p^n)). 
\label{StabilityStr}
\end{equation}

We chose to discretize equation~\eqref{StabilityStr} in time using a $\theta$-scheme discussed in~\cite{glowinski2003finite}:
\begin{equation} 
\rho_s h \frac{\eta^{n+1}-2 \eta^n+\eta^{n-1}}{\triangle t^2} + \mathcal{L} (\theta \eta^{n+1} +(1-2 \theta) \eta^n+\theta \eta^{n-1}) =  
\beta (p_{ext}^{n+1} + \mathcal S (\beta p^n)). 
\label{StabilityStr2}
\end{equation}
It was shown in \cite{glowinski2003finite} that, for a given fixed right hand-side (source term), this scheme is stable for all $0 \le \theta \le 1/2$. 
Thus, we have unconditional stability with respect to the arbitrary ratios of the fluid and structure densities, provided that the 
right hand-side of this equation converges as $n\to\infty$.

A crucial point to observe here is that the right hand-side of this equation, which comes from the pressure loading,
is given by an iterative procedure, and can be written entirely in terms of the initial pressure, the external pressure ($p_{in/out}$),
and the operator $\mathcal S$ whose maximum eigenvalue, as we shall show below, is always less than 1,
for all the choices of $\rho_f$ and $\rho_s h$.  Moreover, we will show below that the right hand-side converges, as the number
of iterations $n\to\infty$,
if $0 \le \beta \le 1$, for all the choices of $\rho_f$ and $\rho_s h$.

To analyze the right hand-side of equation \eqref{StabilityStr}, we first study the eigenvalues of operator $\mathcal S$.
As we shall see below, it is convenient to express the eigenvalues of $\mathcal S$ via the eigenvalues of the
``Neumann to Dirichlet'' operator $\mathcal M_A : H^{-1/2}(\Gamma) \to H^{1/2}(\Gamma)$ 
which is defined to be the trace on $\Gamma$ of the operator $\mathcal R : H^{-1/2}(\Gamma) \to Q$:
\begin{equation}
\mathcal M_A w = \mathcal R w|_{\Gamma},
\label{neumantodir}
\end{equation}
where $\mathcal R$ associates to every $w \in H^{-1/2}(\Gamma)$
the solution  $\mathcal {R}w \in Q$ of the following (pressure) problem:
\begin{equation}
 \left\{\begin{array}{rcll}
 \displaystyle{-\triangle \mathcal {R}w}  &=&0,   \quad  &\textrm{in} \; \Omega, \\ 
  \displaystyle{\mathcal {R}w }  &=& 0  \quad &\textrm{on} \; \Gamma_{in/out}, \\
 \displaystyle{\frac{\partial \mathcal {R}w}{\partial \boldsymbol n}}  &=& 0 \quad &\textrm{on} \; \Gamma_{b}, \\[0.3cm] 
\displaystyle{\frac{ \partial \mathcal {R}w}{\partial \boldsymbol n}} &=& w    \quad &\textrm{on}\; \Gamma.
\end{array} \right.  \nonumber
\label{defR}
\end{equation}
It can be shown (see \cite{causin2005added}) that operator ${\mathcal M}_A : H^{-1/2}(\Gamma) \to H^{1/2}(\Gamma) $ is continuous,
and that ${\mathcal M}_A$ is compact, self-adjoint, and positive on $L^2(\Gamma)$. Moreover, 
the eigenvalues $\mu_i$ of $\mathcal M_A$ are decreasing to zero ($\mu_i = L/(i\pi \tanh (i \pi R/L))$, i = 1,2,...), with the maximum eigenvalue $\mu_{max}$ given by
$$
\mu_{max} = \mu_1 = \frac{L}{\pi \tanh \bigg( \displaystyle{\frac{\pi R}{L}\bigg)}}.
$$
We will use this knowledge to calculate the eigenvalues of operator $\mathcal S$.
Let $\mu$ be an eigenvalue of operator ${\mathcal M}_A$. Then, there exists a vector $v \ne 0$ such that
$$
{\mathcal M}_A v = \mu v.
$$
Recall that, by definition of $\mathcal M_A$, $\mathcal M_A v = \mathcal R v|_{\Gamma}$ and 
$\displaystyle{\frac{ \partial \mathcal {R}v}{\partial \boldsymbol n}} = v$.
Using this, we calculate 
$$
\mathcal R v|_{\Gamma} + \frac{\rho_s h}{\rho_f} \frac{\partial\mathcal R v}{\partial n}|_{\Gamma}= (\mu + \frac{\rho_s h}{\rho_f} ) v |_{\Gamma}.
$$
This implies that $\mathcal R v$ also satisfies the following Robin problem:
\begin{equation}
 \left\{\begin{array}{rcll}
 \displaystyle{-\triangle \mathcal {R}v}  &=&0,   \quad  &\textrm{in} \; \Omega, \\ 
  \displaystyle{\mathcal {R}v}   &=& 0  \quad &\textrm{on} \; \Gamma_{in/out}, \\
 \displaystyle{\frac{\partial \mathcal {R}v}{\partial \boldsymbol n}}  &=& 0 \quad &\textrm{on} \; \Gamma_{b}, \\[0.3cm] 
\displaystyle{\mathcal {R}v + \frac{\rho_s h}{\rho_f} \frac{ \partial \mathcal {R}v}{\partial \boldsymbol n}} &=& \displaystyle{\big(\mu  + \frac{\rho_s h}{\rho_f} \big) v }   \quad &\textrm{on}\; \Gamma.
\end{array} \right.  \nonumber
\end{equation}
We now notice that this is precisely the problem defined by the operator $\mathcal P$ in \eqref{P}, with the data on $\Gamma$ given by $(\mu + \frac{\rho_s h}{\rho_f} ) v$.
Thus:
\begin{equation*}
 \mathcal {R}v =  \mathcal {P} \big(\mu + \frac{\rho_s h}{\rho_f} \big)v,
\end{equation*}
and therefore, the traces on $\Gamma$ satisfy:
\begin{equation*}
 \mathcal {M_A}v = \mathcal {S}  \big(\mu  +  \frac{\rho_s h}{\rho_f} \big) v = \big(\mu  +  \frac{\rho_s h}{\rho_f} \big)\mathcal {S} v.
\end{equation*}
Since ${\mathcal M}_A v = \mu v$ we finally get
$$
\mu v = \big(\mu  +  \frac{\rho_s h}{\rho_f} \big)\mathcal {S} v.
$$
Therefore, $v$ is also an eigenvector for $\mathcal S$, and the corresponding eigenvalue $\lambda$ satisfies:
\begin{equation}\label{eigenval}
\lambda = \frac{\mu}{\mu + \frac{\rho_s h}{\rho_f}}.
\end{equation}
Thus, we have shown that the eigenvalues $\lambda_i$ of $ \mathcal {S}$ 
can be expressed using the eigenvalues $\mu_i$ of $ \mathcal {M_A}$ as
\begin{equation}\label{lambda}
\lambda_i = \frac{\mu_i}{\mu_i+\displaystyle{\frac{\rho_s h}{\rho_f}}}, i = 1,2,....
\end{equation}

We now use this information to study the right hand-side of equation \eqref{StabilityStr2}.
Since $\mathcal S$ is compact, there exists an orthonormal basis of $L^2(\Gamma)$ composed of the eigenvectors $\{z_j\}$ of $\mathcal S$.
We thus expand the solution $\eta$, the external pressure data $p_{ext}$, and $p^0$, in this basis:
$$
\eta^n = \sum_{j} (\eta^n)_j z_j,\quad
p_{ext}^n = \sum_{j} (p_{ext}^n)_j z_j,\quad
p_0 = \sum_{j} (p_0)_j z_j.
$$
Then, from \eqref{StabilityStr2}, for each $j$, the Fourier coefficients satisfy the following equation:
\begin{equation} 
\begin{array}{c}
\displaystyle{\rho_s h \frac{(\eta^{n+1})_j -2 (\eta^n)_j+(\eta^{n-1})_j}{\triangle t^2} + \mathcal{L} (\theta (\eta^{n+1}_r)_j +(1-2 \theta) (\eta^n)_j+\theta (\eta^{n-1})_j)} \\
=  \beta ((p_{ext}^{n+1})_j + \mathcal S (\beta (p^n))_j). 
\end{array}
\label{StabilityStr3}
\end{equation}
The right hand-side of this equation is equal to 
\begin{eqnarray*}
&\beta ((p_{ext}^{n+1})_j + \mathcal S (\beta (p^n))_j) =\\
& \beta  (p_{ext}^{n+1})_j + \sum_{i=1}^n  \beta^{i+1} \lambda_j^i(p_{ext}^{n+1-i})_j + \beta^{n+2} \lambda_j^{n+1} (p_0)_j.
\end{eqnarray*}
As $n\to\infty$, the series that defines the right hand-side converges if 
$$
|\beta \lambda_j| < 1.
$$
From \eqref{lambda} we see that all $\lambda_j$ are strictly less than one, which implies that
the right hand-side converges if
\begin{equation}\label{beta}
0 \le \beta \le 1.
\end{equation}

We have shown the following result:
\begin{theorem}
The  kinematically coupled $\beta$-scheme applied to a class of FSI problems represented by
the simple benchmark problem \eqref{StabilityStokes}, is unconditionally stable 
for each $\beta$ such that $0 \le \beta \le 1$.
\end{theorem}

\section{Numerical results}\label{NumericalResults}
We present two numerical examples which show stability of the kinematically coupled $\beta$-scheme
for a {\bf fully nonlinear FSI problem}. The parameter values in the two examples are well within the
range for which the classical Dirichlet-Neumann schemes are unstable.  
Numerical simulations in both cases are preformed on a benchmark problem  
by Formaggia et al.~\cite{formaggia2001coupling}
used for testing of several FSI algorithms in hemodynamics applications~\cite{badia,nobile2001numerical,badia2008splitting,quaini2009algorithms,guidoboni2009stable}.

Recall that results of \cite{causin2005added} show that the classical loosely coupled, Dirichlet-Neumann scheme is unconditionally unstable if 
\begin{equation}
\frac{\rho_s h}{\rho_f \mu_{max}} < 1,
\label{DNcond}
\end{equation}
where $\mu_{max}$ is the maximum eigenvalue of the added mass operator given by $$\mu_{max}=\frac{L}{\pi \tanh\big( \frac{\pi R}{L}\big)}.$$
For the parameters given in Table~\ref{T1}, the value of $\mu_{max}$ is $7.46$, 
so the critical value for the structure density  is $\rho_s = 74.6$ g/cm$^3$. 
Therefore,  the classical Dirichlet-Neumann scheme is unconditionally unstable if
$\rho_s < 74.6$ g/cm$^3$. 
In Example 1, we take the density of the structure to be 70-times smaller, given by the physiologically relevant
value of  $\rho_s=1.1$ g/cm$^3$. 
We compare our results with the results obtained using a monolithic scheme by Badia, Quaini and Quarteroni \cite{badia2008splitting,quaini2009algorithms} 
showing excellent agreement. 
In Example 2, we choose an even smaller structural density, $\rho_s = 0.55$ g/cm$^3$, 
which is roughly half the value of the fluid density, therefore, capturing the case when 
$\rho_s < \rho_f$. Our results show that the scheme is stable even when the structure is
lighter than the fluid. 

\begin{figure}[ht!]
 \centering{
 \includegraphics[scale=0.65]{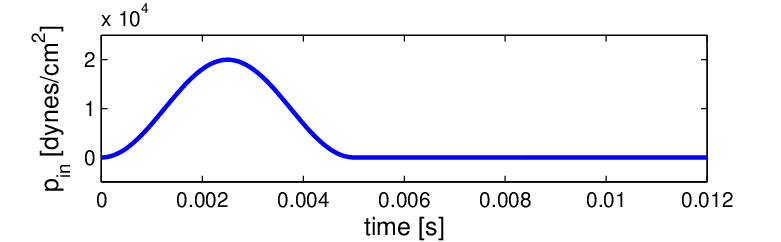}
 }
 \caption{The inlet pressure pulse for Examples 1 and 2. The outlet pressure is kept at 0. }
\label{pressure_pulse}
 \end{figure}

\vskip 0.1in
\noindent
{\bf{The benchmark problem \cite{formaggia2001coupling}:}}\\
The benchmark problem consists of solving the fully nonlinear FSI problem \eqref{NS}-\eqref{Coupling1} 
with the values of the coefficients for the fluid problem given in Table~\ref{T1},
and the values of the structural coefficients given in Table~\ref{Tcoeff}.
{{
\begin{center}
\begin{table}[ht!]
{\small{
\begin{tabular}{|l l | l l |}
\hline
\textbf{Parameters} & \textbf{Values} & \textbf{Parameters} & \textbf{Values}  \\
\hline
\hline
\textbf{Radius} $R$ (cm)  & $0.5$  & \textbf{Length} $L$ (cm) & $6$  \\
\hline
\textbf{Fluid density} $\rho_f$ (g/cm$^3$)& $1$ &\textbf{Dyn. viscosity} $\mu$ (poise) & $0.035$    \\
\hline
\textbf{Young's mod.} $E $(dynes/cm$^2$) & $0.75 \times 10^6$   & \textbf{Wall thickness} $h$ (cm) & $0.1$  \\
\hline
\textbf{Poisson's ratio} $\sigma $ & $0.5$  &  \\
\hline
\end{tabular}
}}
\caption{Parameter values used in Examples 1 and 2 presented in this section.}
\label{T1}
\end{table}
\end{center}
}}
\begin{table}[ht]
\begin{center}
\begin{tabular}{| l  l | l  l |}
\hline
\textbf{Coefficients} & \textbf{Values}  & \textbf{Coefficients} & \textbf{Values} \\
\hline
\hline
\textbf{$C_0$}  & $4 \times 10^5$  &  \textbf{$C_1$} & $2.5 \times 10^4$\\
\hline
\textbf{$D_0$}  & $0$  &  \textbf{$D_1$} & $0.01$\\
\hline
\end{tabular}
\caption{The values of the coefficients in~\eqref{StructureProblem1} used in Examples 1 and 2 below.}
\label{Tcoeff}
\end{center}
\end{table}
The flow is driven by the time-dependent pressure data:
\begin{equation}\label{pressure}
 p_{in}(t) = \left\{\begin{array}{l@{\ } l} 
\frac{p_{max}}{2} \big[ 1-\cos\big( \frac{2 \pi t}{t_{max}}\big)\big] & \textrm{if} \; t \le t_{max}\\
0 & \textrm{if} \; t>t_{max}
 \end{array} \right.,   \quad p_{out}(t) = 0 \;\forall t \in (0, T),
\end{equation}
where $p_{max} = 2 \times 10^4$ (dynes/cm$^2$) and $t_{max} = 0.005$ (s). 
The graph of the inlet pressure data versus time is shown in Figure~\ref{pressure_pulse}.
The inlet and outlet boundary conditions for the structure are the absorbing boundary conditions:
\begin{align}
\frac{\partial \eta_r}{\partial t} - \sqrt{\frac{C_1}{\rho_s}} \frac{\partial \eta_r}{\partial z} & = 0 \quad \textrm{at} \; z=0 \label{abs1} \\ 
\frac{\partial \eta_r}{\partial t} + \sqrt{\frac{C_1}{\rho_s}} \frac{\partial \eta_r}{\partial z} & = 0 \quad \textrm{at} \; z=L. \label{abs2}
\end{align}
To deal with the motion of the fluid domain, we used the ALE mapping defined by the harmonic extension of the mapping that maps the boundary of 
$\Omega$ onto the boundary of $\Omega(t)$, for a given time $t$.

Parameter $\beta$, introduced in \eqref{beta},  which appears in Step 1 and Step 3 of our numerical scheme,
and is independent of time,  was taken to be $\beta=1$.
It was numerically observed in~\cite{MarSun} that the change in $\beta$ is associated with the change in accuracy of the scheme
(not the stability), 
where the value of $\beta = 1$ provides the highest accuracy for this benchmark problem.
We believe that the main reason for the gain 
in accuracy at $\beta = 1$ is the strong coupling between the fluid pressure (which incorporates the leading effect of the fluid loading onto the structure)
and the structure elastodynamics, which is established for $\beta = 1$ in Step 3 of the splitting, described above.

\subsection{Example 1: $\rho_s=1.1$ g/cm$^3$. }
We solved the fully nonlinear benchmark problem described above over the time interval $[0,0.012]$s.
This time interval was chosen to be the same as in  \cite{formaggia2001coupling}. The end-point of this time interval corresponds
to the time it takes a forward-moving pressure wave to reach the end of the fluid domain.
The numerical results obtained using the kinematically coupled $\beta$-scheme were compared with 
the numerical results obtained using the monolithic scheme by Badia, Quaini and Quarteroni~\cite{badia2008splitting,quaini2009algorithms}.
\begin{figure}[ht!]
 \centering{
 \includegraphics[scale=0.65]{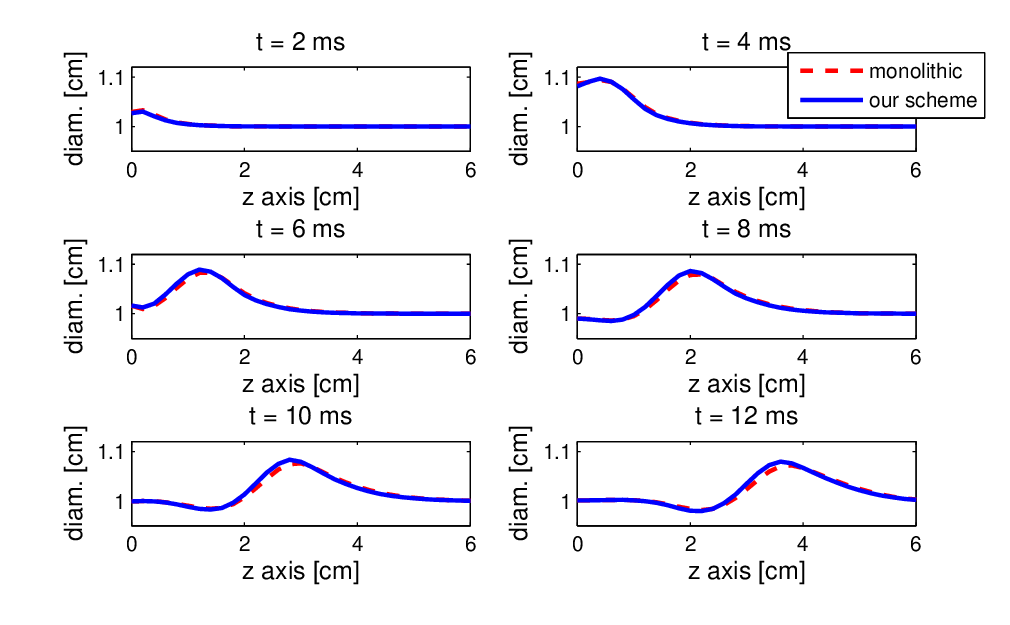}
 }
 \caption{Example 1: Diameter of the tube computed with the kinematically coupled $\beta$-scheme (solid line)
 and with the monolithic scheme used by Quaini in~\cite{badia2008splitting,quaini2009algorithms} (dashed line). The time step $\triangle t = 10^{-4}$ is used in both cases. }
\label{displfig}
 \end{figure}
\begin{figure}[ht!]
 \centering{
 \includegraphics[scale=0.65]{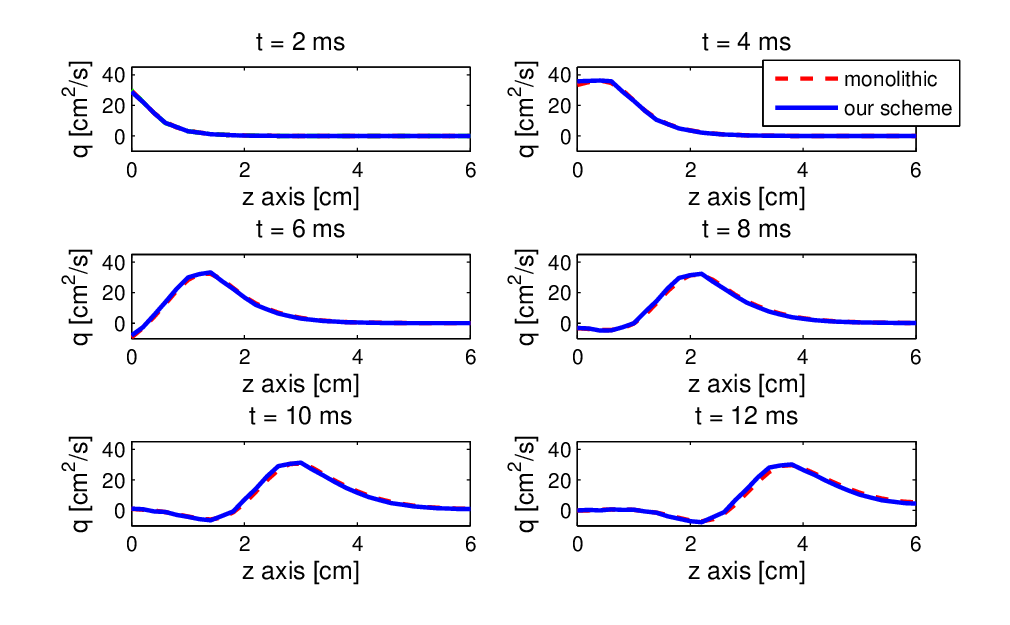}
 }
 \caption{Example 1: Flowrate computed with the kinematically coupled $\beta$-scheme (solid line)
 and  with the monolithic scheme used by Quaini in~\cite{badia2008splitting,quaini2009algorithms} (dashed line).
 The time step $\triangle t = 10^{-4}$ in both cases. }
\label{flow}
 \end{figure}
\begin{figure}[ht!]
 \centering{
 \includegraphics[scale=0.65]{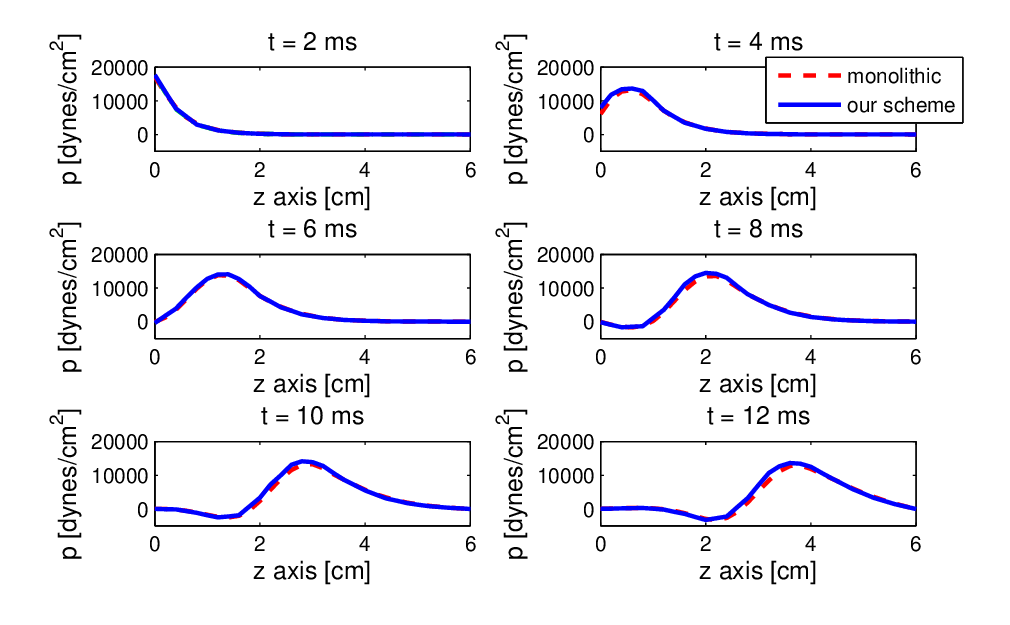}
 }
 \caption{Example 1: Mean pressure computed with the kinematically coupled $\beta$-scheme (solid line)
 and with the monolithic scheme used by Quaini in~\cite{badia2008splitting,quaini2009algorithms}  (dashed line).
 The time step $\triangle t = 10^{-4}$ is used in both cases. }
\label{press}
 \end{figure}
Figures~\ref{displfig},~\ref{flow} and~\ref{press} show a comparison between tube diameter, flowrate and mean pressure, respectively, at six different times. 
These results were obtained on the same mesh as the one used for the monolithic scheme in~\cite{quaini2009algorithms}, containing $31 \times 11$ $ \mathbb{P}_1$ fluid nodes. 
More preciesely, we used an isoparametric version of the Bercovier-Pironneau element spaces, also known as the
$ \mathbb{P}_1$-iso-$ \mathbb{P}_2$ approximation, in which a coarse mesh is used for the pressure (mesh size $h_p$) and 
a fine mesh for the velocity (mesh step $h_v=h_p/2$).

\begin{figure}[ht!]
 \centering{
 \includegraphics[scale=0.65]{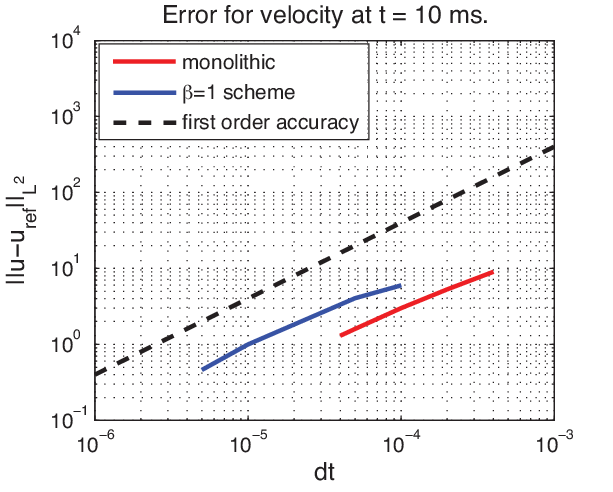}
 \includegraphics[scale=0.65]{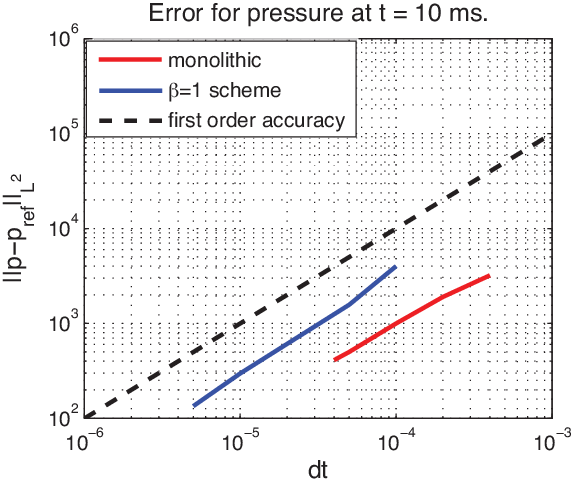} 
 \includegraphics[scale=0.65]{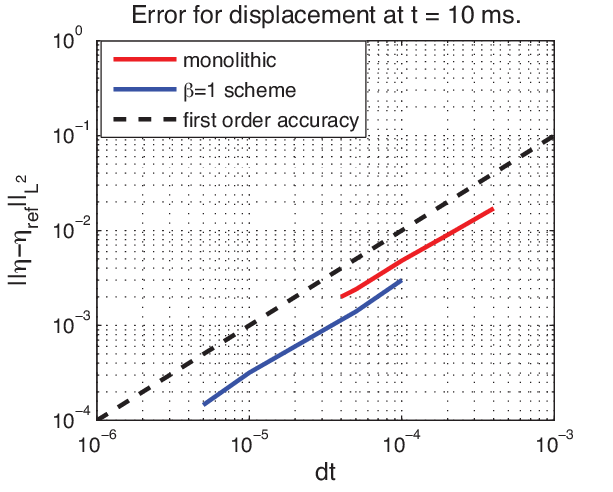}
 }
 \caption{Example 1: The figures show absolute errors compared with the monolithic scheme. 
 Top left: Error for fluid velocity at t=10 ms. Top right: Error for fluid pressure at t=10 ms. Bottom: Error for displacement at t=10 ms,
 \cite{MarSun}. Both schemes are first-order accurate in time, with comparable accuracy (notice the higher accuracy of the kinematically-coupled $\beta$ scheme
 for the displacement).}
\label{error}
\label{ACCURACY}
 \end{figure}

The time step used was $\triangle t = 10^{-4}$, which is the same as the time step used for the monolithic scheme.
Due to the splitting error,
it is well-known that classical splitting schemes usually require smaller time step to achieve the accuracy comparable to 
that of monolithic schemes.
However, the new splitting with $\beta = 1$  allows us to use the same time step as in the monolithic method,
obtaining comparable accuracy, as shown in~\cite{MarSun} and in Figure~\ref{ACCURACY}.

The reference solution was defined to be the one obtained with $\triangle t = 10^{-6}$. 
We calculated the relative and absolute $L^2$ errors for the velocity, pressure and displacement between the reference solution and the solutions obtained using 
$\triangle t = 5\times 10^{-6}, 10^{-5},\; 5 \times 10^{-5}$ and $10^{-4}$.
Table~\ref{T4} shows the relative error and the convergence rates for the pressure, velocity, and displacement obtained by the kinematically-coupled $\beta$-scheme. 
The graphs in Figure~\ref{ACCURACY} show the absolute error and the convergence in time of the kinematically-coupled $\beta$-scheme and the monolithic scheme
by Badia, Quaini and Quarteroni \cite{badia2008splitting,quaini2009algorithms}.
 \begin{table}[ht!]
\begin{center}
{\scriptsize{
\begin{tabular}{| l  c  c  c  c  c  c |}
\hline
$ \triangle t $ & $\frac{||p-p_{ref}||_{L^2}}{||p_{ref}||_{L^2}} $ & $L^2$ order & $\frac{||\boldsymbol u-\boldsymbol u_{ref}||_{L^2}}{||\boldsymbol u_{ref}||_{L^2}}$  &$L^2$ order & $ \frac{||\boldsymbol \eta - \boldsymbol \eta_{ref}||_{L^2}}{||\boldsymbol \eta_{ref}||_{L^2}} $ & $L^2$ order \\
\hline
\hline
$10^{-4}$ & $ 0.0251$  & - & $0.0223$  & - & $0.0392$  & - \\ 
\hline  
$5 \times 10^{-5}$ & $0.013$  & 1.35 & $0.0151$  & 0.56 & $0.0175$  & 1.1 \\ 
\hline   
$10^{-5}$ & $0.0024$  & 1.04 & $0.0038$  & 0.87 & $0.0038$  & 0.92 \\ 
\hline   
$5 \times 10^{-6}$ & $0.0011$  & 1.14 & $0.0017$  & 1.12 & $0.0017$  & 1.13 \\ 
\hline 
\end{tabular}
}}
\end{center}
\caption{Example 1: Convergence in time calculated at $t = 10$ ms.}
\label{T4}
\end{table}

\if 1=0
The kinematically coupled scheme was shown numerically to be first-order accurate  in time and second order accurate in space~\cite{guidoboni2009stable}. 
Second-order accuracy in space is retained in the current kinematically coupled $\beta$-scheme, since the same spatial discretization
of the underlying operators $A_i, i = 1, 2, 3,$ was used in the present manuscript as in \cite{guidoboni2009stable}.
 \begin{table}[ht!]
\begin{center}
{\scriptsize{
\begin{tabular}{| l  c  c  c  c  c  c |}
\hline
$ \triangle t $ & $\frac{||p-p_{ref}||_{L^2}}{||p_{ref}||_{L^2}} $ & $L^2$ order & $\frac{||\boldsymbol u-\boldsymbol u_{ref}||_{L^2}}{||\boldsymbol u_{ref}||_{L^2}}$  &$L^2$ order & $ \frac{||\boldsymbol \eta - \boldsymbol \eta_{ref}||_{L^2}}{||\boldsymbol \eta_{ref}||_{L^2}} $ & $L^2$ order \\
\hline
\hline
$10^{-4}$ & $ 0.0251$  & - & $0.0223$  & - & $0.0392$  & - \\ 
\hline  
$5 \times 10^{-5}$ & $ 1.57 \textrm{e}+03$  & 1.35 & $4.05$  & 0.56 & $0.0014$  & 1.1 \\ 
\hline   
$10^{-5}$ & $ 296.36$  & 1.04 & $ 1.0$  & 0.87 & $3.17 \textrm{e}-04$  & 0.92 \\ 
 & $( 7.27 \textrm{e}+03)$  &(0.95) & $(16.27)$  &(0.97) & $ (0.00576)$  & $(0.95)$\\ 
\hline   
$5 \times 10^{-6}$ & $ 134.33$  & 1.14 & $0.46$  & 1.12 & $1.45 \textrm{e}-04$  & 1.13 \\ 
 & $(3.3 \textrm{e}+03)$  & (1.14) & $ (7.36)$  & $(1.14)$ & $(0.0026)$  & $(1.14)$ \\ 
\hline 
\end{tabular}
}}
\end{center}
\caption{Example 1: Convergence in time calculated at $t = 10$ ms. }
\label{T4}
\end{table}
However, due to the new time-splitting, the accuracy in time has changed.
Indeed, here we show that this is the case by studying time-convergence of our scheme. 
Figure~\ref{error} shows a comparison between the time convergence of our scheme, the kinematically coupled scheme, and the monolithic scheme used in~\cite{quaini2009algorithms}. 
The reference solution was defined to be the one obtained with $\triangle t = 10^{-6}$. 
We calculated the absolute $L^2$ error for the velocity, pressure and displacement between the reference solution and the solutions obtained using 
$\triangle t = 5\times 10^{-6}, 10^{-5},\; 5 \times 10^{-5}$ and $10^{-4}$.
than the error obtained using the classical kinematically coupled scheme, and is comparable to the error obtained by the monolithic scheme. 
The values of the convergence rates for pressure, velocity, and displacement, calculated using the kinematically coupled schemes, are shown in Table~\ref{T4}.
\fi

\subsection{Example 2: $\rho_s = 0.55$ g/cm$^3$.}
\begin{figure}[ht!]
 \centering{
 \includegraphics[scale=0.65]{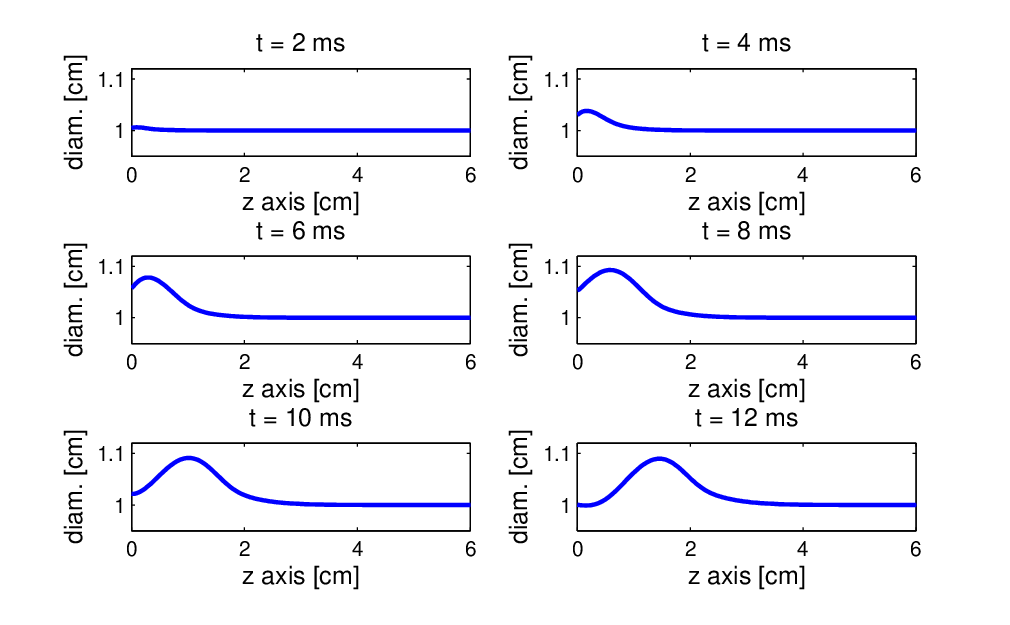}
 }
 \caption{Example 2: Diameter of the tube computed with the kinematically coupled $\beta$-scheme,
 obtained with the time step $\triangle t = 10^{-4}$ . }
\label{displfig2}
 \end{figure}

\begin{figure}[ht!]
 \centering{
 \includegraphics[scale=0.65]{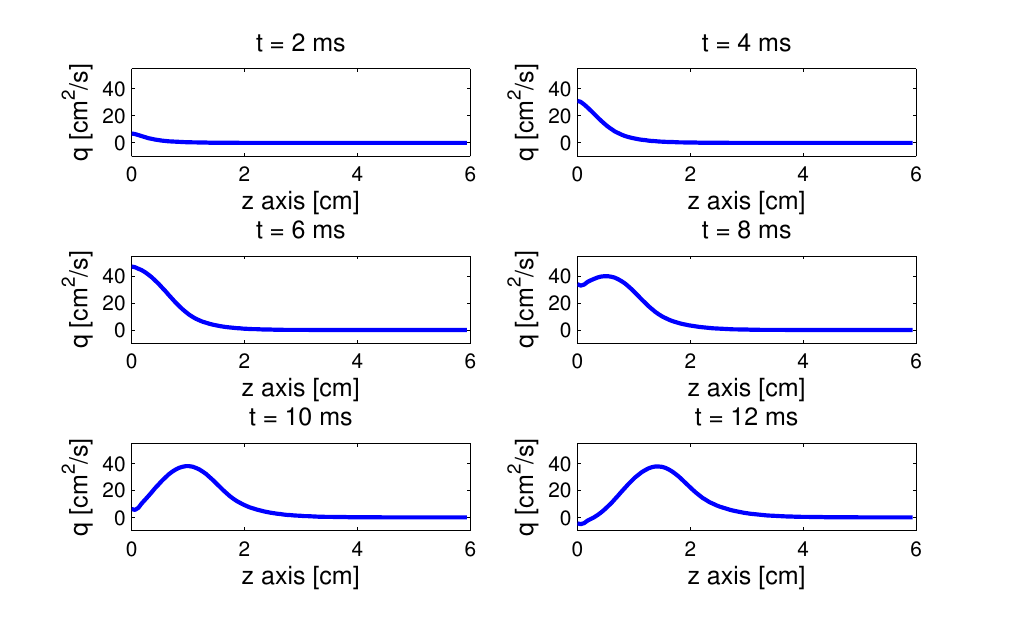}
 }
 \caption{Example 2: Flowrate computed with the kinematically coupled $\beta$-scheme,
 obtained with the time step $\triangle t =  10^{-4}$. }
\label{flow2}
 \end{figure}

\begin{figure}[ht!]
 \centering{
 \includegraphics[scale=0.65]{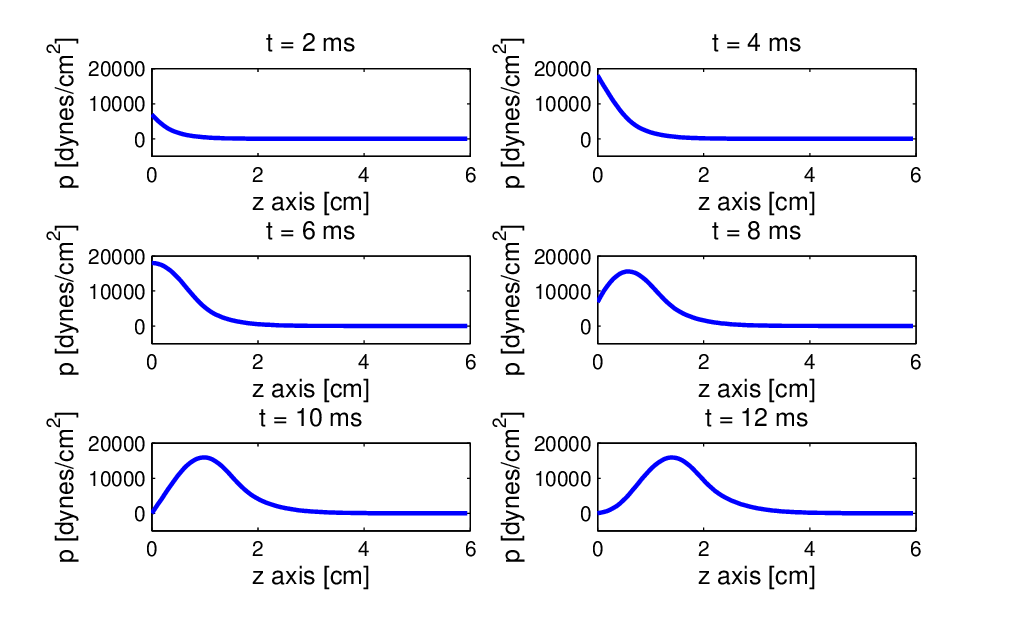}
 }
 \caption{Example 2: Mean pressure computed with the kinematically coupled $\beta$-scheme,
 obtained with the time step $\triangle t = 10^{-4}$. }
\label{press2}
 \end{figure}

We consider the same test case as in Example 1, but now with the value for the structure density $\rho_s = 0.55$ g/cm$^3$. 
The problem is solved on the time interval $[0,0.012]$s and with the same mesh as the one used in Example~1. 
Figures~\ref{displfig2},~\ref{flow2} and~\ref{press2} show the values of the tube diameter, flowrate and mean pressure, respectively, at six different times,
obtained with $\triangle t = 10^{-4}$. 

Finally, to show the convergence of the scheme, we define the reference solution to be the one obtained with $\triangle t = 10^{-6}$. 
We calculated the relative $L^2$ errors for the velocity, pressure and displacement between the reference solution and the solutions obtained using 
$\triangle t = 5\times 10^{-6}, 10^{-5},\; 5 \times 10^{-5}$ and $10^{-4}$.
Table~\ref{T5} shows the convergence rates for the pressure, velocity, and displacement obtained by 
the kinematically-coupled $\beta$-scheme. 
Thus, we see that the scheme performs well for a range of parameters for which the classical loosely-coupled, Dirichlet-Neumann scheme is unconditionally unstable.
 \begin{table}[ht!]
\begin{center}
{\scriptsize{
\begin{tabular}{| l  c  c  c  c  c  c |}
\hline
$ \triangle t $ & $\frac{||p-p_{ref}||_{L^2}}{||p_{ref}||_{L^2}} $ & $L^2$ order & $\frac{||\boldsymbol u-\boldsymbol u_{ref}||_{L^2}}{||\boldsymbol u_{ref}||_{L^2}}$  &$L^2$ order & $ \frac{||\boldsymbol \eta - \boldsymbol \eta_{ref}||_{L^2}}{||\boldsymbol \eta_{ref}||_{L^2}} $ & $L^2$ order \\
\hline
\hline
$10^{-4}$ & $ 0.0239$  & - & $0.0427$  & - & $0.0749$  & - \\ 
\hline  
$5 \times 10^{-5}$ & $ 0.0096$  & 1.32 & $0.0286$  & 0.58 & $0.0408$  & 0.87 \\ 
\hline   
$10^{-5}$ & $0.0017$  & 1.06 & $0.0067$  & 0.9 & $0.0079$  & 1.021 \\ 
\hline   
$5 \times 10^{-6}$ & $7.72 \textrm{e}-04$  & 1.15 & $0.0031$  & 1.13 & $0.0035$  & 1.16 \\ 
\hline 
\end{tabular}
}}
\end{center}
\caption{Example 2: Convergence in time calculated at $t = 10$ ms.}
\label{T5}
\end{table}

\section{Comparison with the Dirichlet-Neumann Scheme}\label{Conclusions}
We conclude this manuscript by comparing the kinematically-coupled $\beta$-scheme, summarized in 
\eqref{StabilityStr}, with the classical
Dirichlet-Neumann scheme exhibiting the ``added mass effect.''
Equation \eqref{StabilityStr} can be solved in an implicit way (written with a slight abuse of notation) as:
\begin{equation} \label{StructureFinal}
\rho_s h \frac{\partial^2 \eta^{n+1}}{\partial t^2} + \mathcal{L} \eta^{n+1} =  \beta (p_{ext}^{n+1} + \mathcal S (\beta p^n)). 
\end{equation}

\begin{figure}[ht] 
 \centering{
 \includegraphics[scale=0.54]{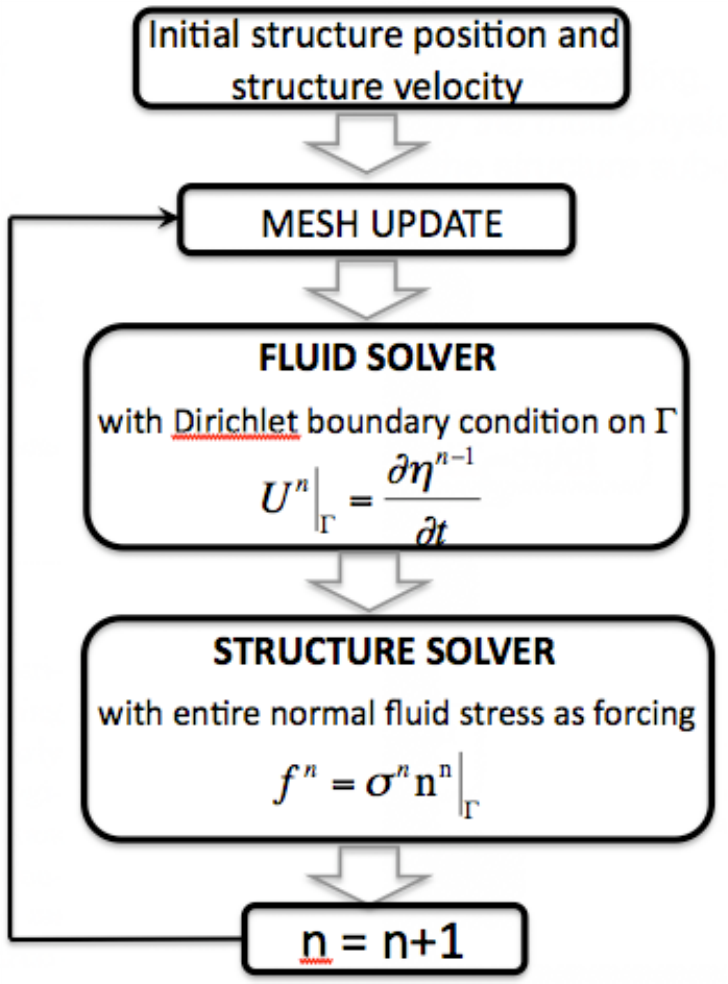}
  \includegraphics[scale=0.4]{kk_beta_scheme.pdf}
 }
 {\bf {\small{ Dirichlet-Neumann Scheme\hskip 0.65in
 Kinematically-Coupled $\beta$-Scheme}}}
 \caption{Left: Block diagram showing the main steps of a Dirichlet-Neumann scheme.
 Right: Block diagram showing the main steps of the kinematically-coupled $\beta$-scheme.}
 \label{classical_scheme}
 \end{figure}

The corresponding equation resulting from the classical Dirichlet-Neumann scheme can be written as follows.
We recall that  the Dirichlet-Neumann scheme solves the FSI problem 
\eqref{FluidProblem}-\eqref{Coupling} by solving the fluid sub-problem \eqref{FluidProblem}
with Dirichlet boundary data for the fluid velocity on $\Gamma$, given in terms of the 
structure velocity $\partial\eta/\partial t$, calculated from the previous time step,
and then uses the fluid stress, calculated in problem \eqref{FluidProblem}, to load
the structure in sub-problem \eqref{StructureProblem}. 
Figure~\ref{classical_scheme} left shows a block diagram 
summarizing the main steps.
 
Using similar ideas as already presented in this manuscript, 
it was shown in \cite{causin2005added} that this kind of partitioned approach
leads to solving the problem for $\eta$ of the form
\begin{equation}\label{st1}
\rho_s h \frac{\partial^2 \eta}{\partial t^2} + \mathcal L \eta = p|_{\Gamma},
\end{equation} 
where
\begin{equation}\label{st2}
p|_{\Gamma} = p_{ext} - \rho_f {\mathcal M}_A  \frac{\partial^2 \eta}{\partial t^2}.
\end{equation}
Here $p_{ext}$ comes from the ``external'' pressure data ($p_{in/out}$), and 
${\mathcal M}_A$ is the ``added mass operator'' introduced in the previous section by formula \eqref{neumantodir}. 

While in the Dirichlet-Neumann loosely coupled scheme the $w$ in the definition of operator ${\cal R}$ (formula \eqref{defR}) corresponds to structure 
inertia $\partial^2 \eta/\partial t^2$,
in our kinematically coupled $\beta$-scheme $w$ corresponds to the fraction of the fluid stress determined by $\beta p$. When $\beta = 0$, which
corresponds to the original kinematically coupled scheme, this term is zero. 
We emphasize again that the structure inertia in the kinematically coupled schemes is taken {\sl implicitly} in the Robin boundary condition,
and not explicitely, as in the Dirichlet-Neumann schemes.

Using the added mass operator, the problem for the structure \eqref{st1} can now be written as
\begin{equation}\label{st3}
(\rho_s h  + \rho_f \mathcal M_A ) \frac{\partial^2 \eta}{\partial t^2} + \mathcal L \eta = p_{ext}.
\end{equation} 
Since in the Dirichlet-Neumann loosely coupled partitioned schemes the pressure is calculated in the fluid sub-problem
using the structure velocity from the previous time step, this implies that, with a slight abuse of notation,
 \eqref{st3} can be written as
 \begin{equation}\label{st4}
\rho_s h  \frac{\partial^2 \eta^{n+1}}{\partial t^2}  + \rho_f \mathcal M_A  \frac{\partial^2 \eta^n}{\partial t^2} +  \mathcal L \eta^n = p_{ext}.
\end{equation} 
More precisely, equation \eqref{st4} means
 \begin{equation}\label{st5}
\rho_s h  \frac{\eta^{n+1} - 2 \eta^n + \eta^{n-1}}{(\Delta t)^2}  +\rho_f \mathcal M_A  \frac{ \eta^n - 2 \eta^{n-1} + \eta^{n-2}}{(\Delta t)^2} +  \mathcal L \eta^n = p_{ext}.
\end{equation} 
It was shown in \cite{causin2005added} that this scheme is unconditionally unstable if 
\begin{equation}
\frac{\rho_s h}{\rho_f \mu_{max}} < 1,
\end{equation}
where $\mu_{max}$ is the maximum eigenvalue of the added mass operator $\mathcal M_A$.
It was also shown in \cite{causin2005added} that the maximum eigenvalue $\mu_{max}$ is associated with the 
aspect ratio of the fluid domain $R/L$. The smaller that aspect ratio (the more slender the domain), the larger the maximum eigenvalue $\mu_{max}$.

Going back to equation \eqref{StructureFinal}, one can now see that there is no added mass operator
associated with the kinematically coupled $\beta$-scheme. The right hand side of equation \eqref{StructureFinal}
is an iterative procedure which can be written in terms of the initial pressure and of the external (in/out) pressure data,
which, we showed, converges for $0 \le \beta \le 1$. 
Thus, for $0 \le \beta \le 1$, the kinematically coupled $\beta$-scheme does not suffer from the added mass effect for
any choice of the parameters in the problem. 
A side-by-side comparison between the two schemes is shown in the block-diagrams  in Figure~\ref{classical_scheme}.

\section{Conclusions}
This manuscript shows the main reasons for the unconditional stability of the kinematically-coupled $\beta$-scheme,
used in \cite{MarSun,CVET,thick_structure,Martina_Biot,guidoboni2009stable,Lukacova,BorSun,BorSunMulti,BorSunStent,BorSun3D} to studycf several different FSI problems in hemodynamics. 
The inclusion of the structure inertia into the fluid sub-problem  implicitly, 
giving rise to a Robin-type boundary condition for the fluid sub-problem, is the main reason 
for the unconditional stability of the scheme. This is in contrast with the classical Dirichlet-Neumann
loosely-coupled schemes, for which the structure inertia is included in the fluid sub-problem
explicitly, via a Dirichlet boundary condition describing the no-slip condition at the fluid-structure interface. 
Modularity, simple implementation, and unconditional stability without the need for sub-iterations at each time-step,
make this scheme particularly appealing for numerical simulation of multi-physics problems arising in biological 
fluid-structure interaction.

\end{document}